\documentclass{comnet}
\bibpunct{[}{]}{,}{n}{,}{;}

\usepackage{natbib}
\usepackage{amsmath} 
\usepackage{amsthm} 
\usepackage{graphics}
\usepackage{algorithmic} 
\usepackage{url} 

\usepackage{graphicx, setspace}
\usepackage{subfig}
\usepackage[table]{xcolor}
\usepackage{float}
\usepackage{multirow}
\usepackage{hhline}

\newcommand{\blue}[1]{\textcolor{blue}{#1}}
\newcommand{\1}{\mathbf{1}}
\newcommand{\esp}{\mathbb{E}}
\newcommand{\argmin}{\mathop{\textrm{argmin}}}

\usepackage[textwidth=3cm]{todonotes}

\begin{document}

\title{Spectral density of random graphs: convergence properties and application in model fitting}

\shorttitle{Random graph model fit via spectral density} 
\shortauthorlist{de Siqueira Santos, Fujita and Matias} 

\author{
\name{Suzana de Siqueira Santos$^*$}
\address{Universidade de S\~ao Paulo,  Instituto de Matem\'atica e Estat\'istica, Departamento de Ci\^encia da Computa\c c\~ ao,  S\~ao Paulo, Brazil}
\address{Funda\c c\~{a}o Getulio Vargas, Escola de Matem\'atica Aplicada, Rio de Janeiro, Brazil\email{$^*$Corresponding author: suzana.santos@fgv.br}}
\name{Andr\'{e} Fujita}
\address{Universidade de S\~ao Paulo,  Instituto de Matem\'atica e Estat\'istica, Departamento de Ci\^encia da Computa\c c\~ ao,  S\~ao Paulo, Brazil}
\and
\name{Catherine Matias}
\address{Sorbonne Universit\'e,  Universit\'e de Paris,
Centre National de la Recherche Scientifique, Laboratoire de Probabilit\'es, Statistique et Mod\'elisation, Paris,
France}}

\maketitle

\begin{abstract}
{Random graph models are used to describe the complex structure of real-world networks in diverse fields of knowledge. Studying their behavior and fitting properties are still critical challenges, that in general, require model-specific techniques. An important line of research is to develop generic methods able to fit and select the best model among a collection. 
Approaches based on spectral density (i.e., distribution of the graph adjacency matrix eigenvalues) appeal to that purpose: they apply to different random graph models. Also, they can benefit from the theoretical background of random matrix theory. 
This work investigates the convergence properties of model fitting procedures based on the graph spectral density and the corresponding cumulative distribution function. We also review the convergence of the spectral density for the most widely used random graph models. Moreover, we explore through simulations the limits of these graph spectral density convergence results, particularly in the case of the block model, where only partial results have been established.}
{random graphs, spectral density, model fitting, model selection, convergence.}
\end{abstract}

\section{Introduction}

The study of real-world networks is fundamental in several areas of knowledge \cite{huckel_quantentheoretische_1931, davidson_gene_2005, pellegrini_protein_2004, van_den_heuvel_exploring_2010}. Analyzing their structure is sometimes challenging, as it may change across time and instances from the same group/phenotype. For example, a functional brain network may change across time and individuals with the same phenotype. Therefore, we can think of a network as a random graph, which is a realization of a random process.

Several models were proposed to describe random processes that generate graphs. For instance, in the Erd\H{o}s-R\'{e}nyi (ER) random graph \cite{erdos_random_1959}, each pair of nodes connects independently at random with probability $p$.
We obtain a generalization when  considering nodes $i$ and $j$ connect independently at random with non-homogeneous probability $p_{ij}$. 
Later, \cite{
frank_cluster_1982, 
holland_stochastic_1983,snijders_estimation_1997} proposed an intermediate setup - block models - 
to generate graphs with non-homogeneous connection probabilities taking a finite number of values. In what follows, we distinguish between the (deterministic) block model (BM) and the stochastic block model (SBM). In the former, pairs of nodes connect independently. In the latter, the connections are only \emph{conditionally} independent (given the nodes' groups).

Models can also characterize different properties of networks. Examples include networks with spatial associations between nodes  \citep[geometric random graph - GRG,][]{penrose_random_2003}, a regular structure ($d$-regular graph - DR),  small-world structure \citep[Watts-Strogatz - WS,][]{watts_collective_1998}, a power-law distribution of the vertex degrees \citep[Barab\'{a}si-Albert - BA,][]{barabasi_emergence_1999}, or occurrence probability depending on pre-defined summary statistics and covariates  \citep[exponential random graph model, ERGM,][]{franck_strauss_1986,chatterjee_estimating_2013}.

Given an empirical network, a significant task is to fit a collection of models through parameter estimation and select the best model in the collection.
There are well-established procedures for estimating the model's parameters \cite{snijders_estimation_1997, ambroise_new_2012, chatterjee_estimating_2013} for some models, such as the ER, SBM, and ERGM. However, those estimators present a lack of generality; they are specific for each model. Every time a new complex network model is proposed, it becomes necessary to develop a specific parameter estimator. Besides, the next question is, among the proposed set of models, which one best fits the observed network? The usual statistical answer relies on model selection criteria in some setups (e.g., ER, SBM) but not all, and are again specific to the model.

In this context, Takahashi et al. proposed fitting and selecting random graph models based on the Kullback-Leibler divergence between graph spectral densities \cite{takahashi_discriminating_2012}. The interest of their approach has been validated on synthetic and real datasets.
The rationale behind this idea is that there is a close relationship between the graph structure and the graph spectrum \cite{von_collatz_spektren_1957}.

The empirical spectral density characterizes some random graph/matrix ensembles. For example, Wigner et al. \cite{wigner_distribution_1958} showed analytically (semicircle law) for ER and DR graphs. In contrast,  Takahashi et al. \cite{takahashi_discriminating_2012} showed empirically for GRG, BA, and WS graphs.

The semicircle law states that a large class of random symmetric matrices' spectral density converges weakly almost surely to a specific distribution called the semicircle distribution \cite{wigner_distribution_1958, grenander_probabilities_1963, arnold_asymptotic_1967}.
We can apply this result to the ER and $d$-regular random graphs under certain conditions \cite{tran_sparse_2013}. However, despite the advances in the study of the convergence of graph spectral densities, we know little about the theoretical properties of the procedures proposed by Takahashi et al.    

In this work, we prove that if we replace the Kullback-Leibler divergence between spectral densities by the $\ell_1$ distance, the parameter estimated by the procedure proposed by Takahashi et al. \cite{takahashi_discriminating_2012} is consistent under the general assumptions:
\begin{enumerate}
    \item the empirical spectral density converges weakly to a limiting distribution,
    \item the map between the parameter space and the limiting distribution is injective and continuous, and
    \item the parameter space is compact (this technical condition could be relaxed).
\end{enumerate}
Furthermore, we show that we can extend the results when considering the eigenvalues' cumulative distribution to fit the models. Note that to obtain theoretical support on the procedures, we need to consider slight modifications of these, e.g., replacing the Kullback-Leibler divergence with the $\ell_1$ distance. This is due to technicalities that may not be overcome yet. 
For example, convergence in total variation norm is necessary for convergence of the Kullback-Leibler divergence (see Corollary 3.1 from \cite{piera_convergence_2009} and Lemma 2.4 from \cite{pinski_kullback--leibler_2015}). This means that theoretical support for the original procedures by \cite{takahashi_discriminating_2012} would require stronger assumptions on the convergence of the empirical spectral density. Our theoretical results (based on the $\ell_1$ distance) are the first to pave the road to explaining the good performance of these procedures.

In Section~\ref{sec:definitions}, we present the notation and main definitions used in this paper. Section~\ref{sec:esd} summarizes the literature on convergence results for the random graph's spectral density under different models. We describe a model-fitting procedure (that slightly differs from \cite{takahashi_discriminating_2012}) in Section~\ref{sec:model_fitting}. We establish its convergence properties in Section~\ref{sec:cv}; see Theorems~\ref{theo:esd} and \ref{theo:cdf}. Then, in Section~\ref{sec:sim}, we evaluate the method performance by simulation experiments. We illustrate applications on real-world data in Section~\ref{sec:real_world_networks}.  Finally, we discuss the results in Section~\ref{sec:discussions} and show our conclusions in Section~\ref{sec:conclusions}. 




%



\section{Notation and definitions}
\label{sec:definitions}
Let $G = (V, E)$ be an undirected graph in which $V = \{1, 2, \dotsc, n\}$ is a set of vertices, and $E$ is a set of edges connecting the elements of $V$. The spectrum of $G$ is the set of eigenvalues of its adjacency matrix $\mathbf{A}=(\mathbf{A}_{ij})_{1\le i,j\le n}$. Since \(G\) is undirected, $\mathbf{A}$ is symmetric, and then its eigenvalues are real values \(\lambda^G_1 \geq \lambda^G_2 \geq \dotsc \geq  \lambda^G_n\).

Let \(\delta_x\) be  the Dirac measure at point $x$. Then, we define the empirical spectral distribution (ESD) as
\begin{equation}
\label{eq:ESD}
\mu^G = \frac{1}{n} \sum_{i=1}^n \delta_{\lambda_i^G/\sqrt{n}}.
\end{equation}
Let $\1\{A\}$ be the indicator function of set $A$. Then, the corresponding empirical eigenvalues' cumulative distribution function (CDF) is
\begin{equation}
    F^G(x) = \frac{1}{n}\sum_{i=1}^n{\1\Big\{\frac {\lambda_i^G} {\sqrt n} \leq x\Big\}}.
\end{equation}

Let $f : \mathbb{R} \mapsto \mathbb{R}$ be a continuous and bounded function. If we  integrate $\mu^{G}$ over $f$ we obtain
\begin{equation}
\label{eq:esd_integration}
\mu^{G}(f) = \int_{\mathbb{R}} f(\lambda) \mu^{G}(d\lambda) = \frac 1 n \sum_{i=1}^n f(\lambda_i^{G}/\sqrt{n}).
\end{equation}
Let $(G_n)_{n\geq 1}$ be a sequence of graphs on the (sequence of) sets of vertices $V$ (in other words, each $G_n$ is a graph over $n$ nodes). We say that $(\mu^{G_n})_{n\geq 1}$ \emph{converges weakly} to $\mu$ (and denote $\mu^{G_n} \Rightarrow \mu$) if for any real bounded and continuous function $f:\mathbb{R} \mapsto \mathbb{R}$, we have $\mu(f) = \lim_{n \rightarrow \infty} \mu^{G_n}(f)$.

In what follows, we will consider a (kernel) function $\phi$ and introduce for any $x\in \mathbb{R}$ the function $\phi_x (\cdot)=\phi(x-\cdot)$. In this case, the quantity in Eq~\eqref{eq:esd_integration} becomes a convolution function:
\begin{equation}
\label{eq:esd_kernel}
x \mapsto \mu^{G_n}(\phi_x) = \int_{\mathbb{R}} \phi(x-\lambda) \mu^{G_n}(d\lambda) = \frac 1 n \sum_{i=1}^n \phi(x-\lambda_i^{G_n}/\sqrt{n}).
\end{equation}


A \emph{random graph} is a random variable (r.v.) taking its values in the set of graphs. A \emph{random graph model} is a collection of graphs that are either finite or countable, together with a probability distribution $P$ on this collection. Let $(\mathcal{G}_n)_{n \geq 1}$ be a sequence of random graphs on the set of vertices \(V = \{1, 2, \dotsc, n\}\). 
Then the eigenvalues of $\mathcal{G}_n$ are random variables, and the ESD $\mu^{\mathcal{G}_n}$ is a random probability measure. Now, there are many ways in which the random measure $\mu^{\mathcal{G}_n}$ may converge weakly. We say that $\mu^{\mathcal{G}_n}$ converges weakly in expectation to a probability measure $\mu$ if for any real bounded and continuous function $f:\mathbb{R} \mapsto \mathbb{R}$, we have
\[
\esp \mu^{\mathcal{G}_n}(f) \mathop{\to}_{n \rightarrow \infty} \esp\mu(f).
\]
Similarly, we say that $\mu^{\mathcal{G}_n}$ converges weakly in $P$-probability to a probability measure $\mu$ if for any real bounded and continuous function $f:\mathbb{R} \mapsto \mathbb{R}$ and any $\epsilon>0$, we have that,
\[
P\left(|\mu^{\mathcal{G}_n}(f) - \mu(f)| > \varepsilon \right) \to 0,
\]
as $n \to \infty$. Finally, if 
\[
P\left(\lim_{n\to \infty}\mu^{\mathcal{G}_n}(f)= \mu(f) \right) = 1
\]
for any bounded and continuous $f:\mathbb{R} \mapsto \mathbb{R}$, we say that $\mu^{\mathcal{G}_n}$ converges weakly $P$-almost surely to $\mu$.

Let $f$ and $g$ be functions of $n$. Then, as $n \to \infty$, we say that $f = O(g)$ if $|f|/|g|$ is bounded from above; $f = o(g)$ if $f/g \to 0$; and $f = \omega(g)$ if $|f|/|g| \to \infty$.

Finally, we introduce the Stieltjes transform~\cite{bai_methodologies_1999}, which is an important tool for studying the limiting spectral distribution of random matrices. Given a probability distribution $\mu$, its Stieltjes transform $s_{\mu}$ is a function on the upper-half complex plane defined for $z \in \mathbb{C}^+$ as 
\[s_\mu(z) = \int_{\mathbb{R}}\frac{d \mu(x)}{x - z}.\]

\section{Random graph models' spectral density}
\label{sec:esd}

In this section, we summarize results on the spectral density of random graphs. First, we describe Wigner's law (also called semicircle law) for random symmetric matrices. We show examples of random graph models for which the law is valid. Then we summarize analytical and empirical results on the spectral density under different models when Wigner's law does not hold. 

\subsection{Wigner's law}

Many results on the limiting ESD of random graphs rely on Wigner's law for random symmetric matrices. Wigner \cite{wigner_distribution_1958} proved that for symmetric matrices whose entries are real-valued independent and identically distributed (i.i.d.) random variables, with mean zero and unit variance, the ESD converges weakly in expectation to 
the semicircle law
$\mu_{sc}$ defined as 
\begin{equation}
\label{eq:semicircule}
\mu_{sc}(dx) = \frac 1 {2 \pi} \sqrt{4-x^2} \1\{ |x| 
\leq 2\}dx.
\end{equation}
Later, Grenander \cite{grenander_probabilities_1963} proved that this convergence holds weakly in probability. Finally, Arnold \cite{arnold_asymptotic_1967} proved that (under the same assumptions) 
the ESD converges weakly almost surely to the semicircle distribution.


\subsection{Erd\H{o}s-R{\'e}nyi random graph model (ER)}
\label{sec:er}

One of the simplest examples of random graph models in terms of construction is the model proposed by \cite{erdos_random_1959}. Given a probability $p$, and the number of vertices $n$, each pair of vertices is connected with probability $p$, independently of the other pairs.

Let $(\mathcal{G}_{p,n})_{n \geq 1}$ be a sequence of ER random graphs, and $\lambda_1 \geq \lambda_2 \geq \cdots \geq \lambda_n$ denote the eigenvalues of the adjacency matrix $\mathbf{A}$ (here, we do not stress that $\mathbf{A}=\mathbf{A}_n$ depends on $n$). Based on  Wigner's semicircle law, it can be proved that if $p = \omega(\frac{1}{n})$, then the ESD of the scaled matrix
$\mathbf{A}/\sqrt{np(1-p)}$, given by 
\begin{equation}
\label{eq:esd_er}
\tilde\mu^{\mathcal{G}_{p,n}} = \frac{1}{n} \sum_{i=1}^{n} \delta_{ \lambda_i/\sqrt{np(1-p)}},
\end{equation}
converges weakly almost surely (and also in probability and expectation) to the 
semicircle
distribution~\citep{tran_sparse_2013}.
Another formulation says the ER graph $\mathcal{G}_{p,n}$ ESD, denoted by $\mu^{\mathcal{G}_{p,n}}$ and defined as in Equation \eqref{eq:ESD}, converges to
\begin{equation}
\label{eq:semicircle_p}
\mu_{p, sc}(dx) = \frac 1 {2 \pi p(1-p)} \sqrt{4p(1-p)-x^2} \1\{ |x| \leq 2\sqrt{p(1-p)}\}dx.
\end{equation}
Note that if $np =O(1)$, the convergence no longer holds. 

\subsection{(Deterministic) Block model (BM)}
The (deterministic) block model was introduced in the 70's by \cite{Breiger_75,White_76}; see also \cite{Faust_Wasserman_81}. As already mentioned, a natural generalization of the ER model is to consider that any pair of nodes $i,j$ connect independently but with its probability $p_{ij}$. From a statistical point of view, such a model has too many parameters. We can interpret the block model as a particular case where we restrict these $p_{ij}$ to a finite number of values. 
More precisely, each node belongs to one of the $M$ groups, determining its probability of connection.  That is, pairs of nodes connect with probabilities $p_{m,l}$, which depend on the respective groups $m,l$ of these nodes.

In what follows, we consider, more specifically, the model studied in \cite{avrachenkov_spectral_2015}. 
Let $\mathcal{G}$ be a random graph with $n$ nodes $V = \{1, 2, \dotsc, n\}$ belonging to $M$ groups $\{\Omega_1, \Omega_2, \dotsc, \Omega_M\}$ with equal size $K=n/M$. Assume that the conditional probabilities of
connections are as follows
\[
  P(\mathbf{A}_{ij}=1 | i \in m, j \in l) = \left\{
    \begin{tabular}[]{cc}
     $ p_0$ & \text{ if } $m\neq l$ ,\\
      $p_m$ & \text{ if } $m = l$ .
    \end{tabular}
    \right.
  \]
Let $\mathbf{A}$ denote the adjacency matrix of $\mathcal{G}$ and $\mathbf{\bar{A}}$ denote its expectation matrix. Let $i \in \Omega_m$ and $j \in \Omega_l$  be two vertices of $\mathcal{G}$, then the $n \times n$ matrix $\mathbf{\bar{A}}$ is such that $\mathbf{\bar A}_{ij} = p_m$, if $m = l$, and  $\mathbf{\bar A}_{ij} = p_0$, otherwise. Denote $p^\star=\max_{1\le m \le M} p_m $, and $\gamma(n) = 1/\sqrt{np^\star(1-p^\star)}$. The centered and normalized adjacency matrix of $\mathcal{G}$ is defined as 
\begin{equation}
\label{eq:Asbm}
\mathbf{\tilde{A}} = \frac{(\mathbf{A} - \mathbf{\bar A})}{\gamma(n)}.
\end{equation}
An expression for the Stieltjes transform 
of the spectral density of $\mathbf{\tilde{A}}$ was obtained by \cite{avrachenkov_spectral_2015}, as described in Proposition~\ref{prop:avra}.

\begin{proposition} 
\label{prop:avra}
\textsc{(Corollary 3 from \citep{avrachenkov_spectral_2015}).} If $\lim_{n \to +\infty} n p_m = +\infty$ and $p_m/p_0 \leq c$ for some constant $c > 0$ and for all $m = 1, 2, \dotsc, M$, then, almost surely, the  ESD associated with $\mathbf{\tilde{A}}$ converges weakly to a distribution function whose Stieltjes transform is
\[s(z) = \sum_{m=1}^M c_m(z),\]
where the functions $(c_m)_{1\le m \le M}$ are the unique solutions to the system of equations
\[c_m(z) = \frac{-1/M}{z + \zeta_m c_m(z) + \zeta_0 \sum_{l\neq m}c_l(z)}, \quad 1\le m \le M,\]
with
\[\zeta_l = \lim_{n \to +\infty}{\frac{p_l(1-p_l)}{p^\star (1-p^\star)}}, \quad 0\le l \le M,\]
that satisfy the conditions
\[\Im(c_m(z))\Im(z) >0, \text{ for } \Im(z) > 0, \quad 1\le m \le M,\]
where $\Im(z)$ denotes the imaginary part of $z$.
\end{proposition}

The recent paper \cite{Zhu_18} proposes explicit moment  formula for the limiting distribution. 
Notice that the results from \cite{avrachenkov_spectral_2015} were obtained only for the centered matrix $\mathbf{\Tilde{A}}$. For the Erd\H{o}s-R\'enyi model, Tran et al. \cite{tran_sparse_2013} proved that the centered and non-centered adjacency matrices respective spectral densities are  approximately identical. To investigate whether the BM graphs have a similar property, we did simulation experiments comparing the ESD of $\mathbf{\tilde{A}}$ and $\mathbf{{A}}/\gamma(n)$. 
We did simulations to obtain the ESD under various scenarios. We show the results in Section~\ref{sec:sim_block_model}.

\subsection{Stochastic block model (SBM)}
\label{sec:sbm}

Similarly to the previous model, the \emph{stochastic block model (SBM)} produces graphs with groups of vertices connected with a particular edge density \cite{fienberg_categorical_1981,holland_stochastic_1983,snijders_estimation_1997}. The main difference between SBM and the deterministic block model lies in a random assignation of each node to one of the groups with some pre-specified group probability. This results in groups with random sizes, while the latter are fixed in BM. Moreover, the entries of an adjacency matrix in BM are independent random variables, while for SBM, these are conditionally independent random variables. Handling the non-independence of the entries of the adjacency matrix in SBM is more challenging. 
Note that many authors abusively call SBM what is only a BM. 
More precisely, 
given a set of vertices $V=\{1,\dots, n\}$, each vertex first picks a group at random among $M$ possibilities independently of the others and with probabilities $(\pi_1,\pi_2, \dots, \pi_M)$. Then conditional on these latent (unobserved) groups, two vertices from groups $m,l$ connect independently of the others with probability $p_{ml}$. 
To our knowledge, no theoretical results about the convergence of the ESD have been obtained for the SBM yet.

\subsection{Configuration random graph model}

The configuration model is a particular instance of the inhomogeneous ER model and has at least two different usages in the literature. Let
$\underline{d}=(d_1,\dots, d_n)$ be a degree sequence of a graph.

The fixed-degree model $FD(\underline{d})$ is the collection of all graphs with this degree sequence and uniform probability on this collection. We obtain this model by simulation: start with a real graph from degree
sequence $\underline{d}$ and use a rewiring algorithm to shuffle the edges. 

The random-degree model $RD(\underline{d})$ considers all graphs over  $n$  nodes such that the entries 
  $\mathbf{A}_{ij}$  are independent with (nonidentical) Bernoulli distribution  $\mathcal{B}(p_{ij})$  where 
  $p_{ij}=   d_i  d_j/C$ and  $C>0$ is a normalizing constant chosen such that   $0\le p_{ij}\le 1$ (for instance, take $C= \max_{i\neq j} d_id_j$).
It is a random graph model where graphs have degrees only approximately given
by $\underline{d}$. Indeed, if $D_i$ denotes the random degree of node $i$, then 
\[
\esp(D_i)= \sum_{j\neq i}\esp(\mathbf{A}_{ij})= \sum_{j\neq i} p_{ij}= \frac {d_i}{C}
\sum_{j\neq i} d_j =\frac {d_i(2|E|-d_i)}{C}. 
\]
Whenever  $d_i$  is not too large and   $C\simeq  2|E|$  we get $\esp(D_i)\simeq d_i$.

Let $\mathbf{A'}$ denote the normalized adjacency matrix:
\[\mathbf{A'} =  \frac{\mathbf{A}}{\sqrt{n}\sum_{i=1}^n{d_i}}.\]
Preciado and Rahimian \cite{preciado_moment-based_2015} proved that for the random-degree model proposed by Chung and Lu \cite{chung_connected_2002}, under specific conditions regarding the expected degree sequence (sparsity, finite moments, and controlled growth of degrees), the ESD of the normalized adjacency matrix $\mathbf{A'}$ converges weakly almost surely to a deterministic density function. Although the analytical expression of this function is unknown, the authors show its moments explicitly.

\subsection{$d$-regular random graph (DR)}
A $d$-regular random graph is a random graph over $n$ nodes with all degrees equal to $d$. In this case, the parameter $d$ is then directly determined by the observed graph. Although there is no need to estimate $d$, studying this random graph's spectral density can help obtain a measure of goodness of fit. This can be useful for graphs that ``are almost'' regular, in which degrees vary due to measurement errors. 

If $d\ge 3$ is fixed with $n$, and $(\mathcal{G}_{d,n})_{n \geq 1}$ is a sequence of $d$-regular random graphs, then the ESD defined as
\begin{equation}
\label{eq:esd_nonscaled}
\tilde \mu^{\mathcal{G}_{d,n}} = \frac{1}{n} \sum_{i=1}^{n} \delta_{ \lambda_i}
\end{equation}
converges weakly almost surely to 
\begin{equation}
\label{eq:esd_regular_dfixed}
\mu_{d}(dx) = \frac{d \sqrt{4(d-1)-x^2}} {2\pi (d^2 -x^2)} \1\{ |x| \le 2 \sqrt{d-1}\}dx.
\end{equation}
Note that $\tilde \mu^{\mathcal{G}_{d,n}}$ is the distribution of the eigenvalues of $\mathbf{A}$, without scaling by $\sqrt{n}$. 
The result dates back to the work \cite{mckay_expected_1981}. It is also valid when considering that $d$ grows slowly with $n$ \citep{dumitriu_sparse_2012,tran_sparse_2013}. The distribution $\mu_d$ is known as a
Kesten-McKay distribution. 
The result may also be rephrased in the following way \citep{bauerschmidt_local_2017}: the ESD of the matrix $\mathbf{A}/\sqrt{d-1}$ converges weakly almost surely to
\begin{equation}
\label{eq:esd_regular}
\tilde \mu_{d}(dx) =\left( 1 + \frac 1 {d-1} - \frac {x^2} {d} \right)^{-1} \frac{\sqrt{4-x^2}}{2\pi} \1\{ |x| \le 2\} dx.
\end{equation}
Dumitriu and Pal \cite{dumitriu_sparse_2012} showed that the ESD of $\mathbf{A}/\sqrt{d-1}$ converges in probability to $\mu_{sc}$ whenever $d$ tends to infinity slowly, namely $d = n^{\epsilon_n}$ with $\epsilon_n=o(1)$.

Let $\mathbf{J}$ be a $n \times n$ matrix containing only ones, and set
\begin{equation}
\label{eq:M}
\mathbf{M} = \left[\frac{d}{n}\left(1-\frac{1}{n}\right)\right]^{-1/2}\left(\mathbf{A} - \frac{d}{n}\mathbf{J}\right).
\end{equation}
 Tran et al. \cite[Theorem 1.5 in][]{tran_sparse_2013} showed that the ESD of $\mathbf{M}/\sqrt{n}$  converges weakly to the semicircle distribution $\mu_{sc}$ when $d$ tends to infinity with $n$.

\subsection{Geometric random graph model}
To construct a geometric random graph (GRG), we randomly draw $n$ vertices on a $D$-dimensional unit cube (or torus) and connect pairs of vertices whose distance is at most $r_n$ (with $r_n \to 0$). A statistical motivation for this model is that the nodes' coordinates represent (latent) attributes. The metric captures the similarity between individuals and translates into connectivity in the graph.

The following results concern the ESD of (unscaled) adjacency matrix in GRG  under different regimes.
Bordenave \cite{Bordenave_08} proved its convergence  when $nr_n^D \to c \in (0, \infty)$, which corresponds to the thermodynamic regime. He also  proposed an approximation of the limit when the constant $c$ becomes large.
Later, Hamidouche et al. \cite{Hamidouche_etal_19} proved that in the so-called connectivity regime, when $\left(\frac{n}{\log n}\right)
r_n^D \to c \in (0, \infty)$ as $n \to \infty$, the ESD converges to the ESD of a deterministic geometric graph. 
Finally, Bordenave \cite{Bordenave_08} also characterizes the ESD of the adjacency matrix, when normalized  by $n$, in the dense regime, i.e, when $r_n^D$ scales as a constant (namely $r_n^D=O(1)$ and $r_n^D=\omega(1)$).

\subsection{Preferential attachment models}
Preferential attachment random graph models dynamically generate graphs by successively adding a new node in a network together with a set of edges. Probability of connecting the new vertex to an existing node is proportional to its current degree. 
One of the best known preferential attachment model is the Barab\'asi-Albert (BA) one~\citep{barabasi_emergence_1999}. In this model, the probability of connecting a vertex to other existing nodes is proportional to their degrees to power a scaling exponent $p_s$.

There is no theoretical result about the convergence of the empirical spectral density for the Barab\'{a}si-Albert model. However, \cite{farkas_spectra_2001} has described the shape of its distribution through simulations. 

In the linear preferential attachment tree model \cite{bhamidi_spectra_2012}, the probability of connecting a new vertex is proportional to the out-degree plus $1+a$, where $a$ is a process parameter.  Under certain conditions, the properly rescaled ESD converges to a deterministic probability function as $n \to +\infty$ \cite{bhamidi_spectra_2012}.

\subsection{Watts-Strogatz random graph model}

Watts-Strogatz (WS) random graphs \cite{watts_collective_1998} present small-world properties, such as short average path lengths (i.e. the majority of vertices can be reached from all other vertices by a small number of steps) and a higher clustering coefficient (i.e. the number of triangles in the graph) than ER random graphs. The process generating WS graphs starts with a regular lattice of size $n$, with each vertex connected to the $K$/2 nearest vertices on each ring's side. For each vertex $i$, and each edge connected to $i$, with probability $p_r$, it replaces the edge by a new one connecting $i$ to a vertex chosen at random.  To the best of our knowledge, all results regarding the limiting ESD of WS graphs are empirical  \cite{farkas_spectra_2001}.

\subsection{Miscellenaous}
Our list of random graphs models is not exhaustive and focuses on the most common ones.  We conclude it by mentioning a very recent result \cite{avrachenkov_2021} on the soft geometric block model. This model combines the block model approach and geometric random graphs. In this reference, the authors describe the limiting ESD in terms of the Fourier transform \citep[see Theorem 1 in][]{avrachenkov_2021}.

\section{Fitting a model}
\label{sec:model_fitting}
Consider a parameterized random graph model $\{ P_\theta ; \theta\in \Theta\}$ where $P_\theta$ is a probability distribution on a collection of graphs, described through the parameter $\theta \in \Theta$, $\Theta\subset \mathbb{R}^d$. 
Let $\mathcal{G}$ be a random graph from some distribution $P$ (that might belong or not to the parametric model $\{ P_\theta ; \theta\in \Theta\}$). We want to estimate a parameter $\theta$ such that $P_\theta$ is close to $P$.
In many random graphs models, the ESD of graphs generated under the distribution $P_\theta$ converges weakly in some sense (for $P_\theta$) to some distribution that we denote $\mu_\theta$ whenever this limit exists.

ESDs are discrete distributions. In general, the limits $\mu_\theta$ are absolutely continuous measures (that is, they have densities with respect to the Lebesgue measure on the spectrum set $\Lambda$). It is thus natural to consider a kernel estimator of the spectrum's density as defined in Eq~\eqref{eq:esd_kernel}. 
We denote by  $\phi_{x,\sigma} (\cdot)=\sigma^{-1}\phi(\frac{x-\cdot} \sigma)$ a kernel function (a function such that $\int \phi(x)dx= 1$) with parameters $x\in \mathbb{R}$ and $\sigma>0$. The kernel estimator of the spectral density with bandwidth $\sigma$ is defined as 
\[
\mu^{\mathcal{G}}(\phi_{x,\sigma}) = \frac 1 {n\sigma} \sum_{i=1}^n \phi\left( \frac{x-\lambda_i^\mathcal{G}/\sqrt{n}} {\sigma}\right) . 
\]
Suppose the bandwidth $\sigma=\sigma_n$  goes to 0 as $n$ tends to infinity. In that case, this estimator should approximate the limiting spectral density $\mu_\theta$ of the graph. 
Note that we choose to focus on the ESD derived from the normalized adjacency matrix, with normalization (or scaling) $\sqrt{n}$, while the literature review from Section~\ref{sec:esd} contains results with different normalization.  This point is further investigated in Section~\ref{sec:sim}.

Now, given a divergence measure $D$ on the set of probability measures over the spectrum of the graphs, 
we estimate $\theta$ by solving the minimization problem

\begin{equation}
\label{eq:parameter_estimator}
\hat\theta = \argmin_{\theta \in \Theta} D(\mu^{\mathcal{G}}(\phi_{\cdot,\sigma}), \mu_\theta(\phi_{\cdot,\sigma})).
\end{equation}

We either use an analytical expression for $\mu_\theta$ 
(when available) or rely on Monte Carlo estimates. In the latter case, we sample $M$ graphs $\mathcal{G}_1, \mathcal{G}_2, \cdots, \mathcal{G}_M$ from the distribution $P_\theta$. Then, for each graph $\mathcal{G}_m$, we numerically obtain $\mu^{\mathcal{G}_m}(\phi_{x,\sigma})$  and we estimate  $\mu_\theta(\phi_{x,\sigma})$ by
\[\hat \mu_\theta(\phi_{x,\sigma}) =  \frac 1 M \sum_{m=1}^M \mu^{\mathcal{G}_m}(\phi_{x,\sigma}) = \frac 1 {nM\sigma} \sum_{m=1}^M \sum_{i=1}^n \phi\left( \frac{x-\lambda_i^{\mathcal{G}_m}/\sqrt{n}} {\sigma}\right).\]

The parameter estimator problem~\eqref{eq:parameter_estimator} is based on the assumption that the random graph $\mathcal{G}$ is generated by a known family of random graphs  $\{ P_\theta; \theta\in \Theta\}$ (the model). That is, we assume that $P \in \{ P_\theta; \theta\in \Theta\}$.
However, usually, the model $\{ P_\theta; \theta\in \Theta\}$ is unknown. Then given a list of random graph models, we choose the one with the best fit according to some objective criterion. If we consider random graphs models with the same complexity (number of parameters), we choose the model with the smallest divergence
($D$) between the spectral densities, as described in Algorithm~\ref{alg:model_selection}. Otherwise, we penalize the divergence $D$ by a function of the number of parameters to prevent from overfitting.

A previous work by \cite{takahashi_discriminating_2012} uses AIC as a criterion to choose a model. Their approach for problem~\eqref{eq:parameter_estimator} is based on the Kullback-Leibler divergence $KL(\mu^{\mathcal{G}}(\phi_x), \mu_\theta(\phi_x))$, where  $\phi$ is the Gaussian kernel. 
In this case, we obtain 
\begin{equation*}
D(\mu^{\mathcal{G}}(\phi_{x,\sigma}), \mu_\theta(\phi_{x,\sigma}))= \frac 1 {n \sigma} \sum_{i=1}^n \int \phi(x-\lambda_i^{\mathcal{G}}/\sqrt{n})\log \left( \frac{ \sum_{i=1}^n  \phi(x-\lambda_i^{\mathcal{G}}/\sqrt{n})} {n \sigma \mu_\theta(\phi_{x,\sigma})} \right) dx.
\end{equation*}

\begin{algorithm}
Procedure for fitting and selecting a random graph model.
\begin{algorithmic} 
\STATE \textbf{Input:} Graph $\mathcal G$, a list of random graph models $\{ P^i_{\theta} ; \theta\in \Theta^i \subset \mathbb{R}^{d}\}$, a finite subset $\tilde\Theta^i \subset \Theta^i$, for $i = 1, 2, \ldots, N$.\\
\STATE \textbf{Output:} Return model and parameter with best fit.
\STATE Min $\leftarrow +\infty$ 
\STATE Argmin $\leftarrow (0, 0)$
\STATE Compute the kernel density estimator $\mu^{\mathcal G}(\phi_{\cdot,\sigma})$.
\FOR{each parameterized random graph model $\{ P^i_{\theta} ; \theta\in \Theta^i\}$, $i = 1, 2, \cdots, N$,}
    \FOR{each $\theta^{j} \in \tilde\Theta^i$}
        \IF{the limiting ESD from $P^i_{\theta^j}$, denoted by $\mu^i_{\theta^{j}}$, is known analytically}
            \STATE $D_{i,j} \leftarrow D(\mu^{\mathcal G}(\phi_{\cdot,\sigma}) , \mu^i_{\theta^{j}}(\phi_{\cdot,\sigma}) )$. 
        \ELSE
            \STATE Sample $M$ graphs $\mathcal{G}_1, \mathcal{G}_2, \cdots, \mathcal{G}_M$ from $P^i_{\theta^j}$. 
            \FOR{each graph $\mathcal{G}_m$} 
                \STATE Compute the kernel density estimator   $\mu^{\mathcal{G}_m}(\phi_{\cdot,\sigma})$.
                \ENDFOR
            \STATE $\hat\mu^i_{\theta^{j}}(\phi_{\cdot,\sigma})  \leftarrow \frac 1 M\sum_{m=1}^M  \mu^{\mathcal{G}_m}(\phi_{\cdot,\sigma})$.
            \STATE $D_{i,j} \leftarrow D(\mu^\mathcal{}(\phi_{\cdot,\sigma}) , \hat\mu^i_{\theta^{j}}(\phi_{\cdot,\sigma}) )$
            \ENDIF
        \IF{$D_{i,j} <$  Min} 
            \STATE Argmin $\leftarrow (i, j)$
            \STATE Min $\leftarrow D_{i,j}$
            \ENDIF
        \ENDFOR
    \ENDFOR
\RETURN Argmin.
\end{algorithmic}
\label{alg:model_selection}
\end{algorithm}


Now we focus on $D$ being the $\ell_1$ distance so that we have 
\begin{equation*}
D(\mu^{\mathcal{G}}(\phi_{\cdot,\sigma}), \mu_\theta(\phi_{\cdot,\sigma}))= \|\mu^{\mathcal{G}}(\phi_{\cdot,\sigma}) - \mu_\theta(\phi_{\cdot,\sigma})\|_1 = \int | \mu^{\mathcal{G}}(\phi_{x,\sigma})-\mu_\theta(\phi_{x,\sigma})| dx .
\end{equation*}
We also suppose that the ESD from the parameterized random graph model is known analytically. 
In that case, the convergence properties of $\hat \theta$ are described in the next section.

\section{Convergence properties}
\label{sec:cv}
In what follows, we consider as divergence the $\ell_1$ distance and denote it  $||\cdot||_1$. We also denote $\mathcal{M}_1(\Lambda)$ the set of probability measures over the spectrum $\Lambda$. 

\begin{theorem}
\label{theo:esd}
Let $\{ P_\theta ; \theta\in \Theta\}$ denote a parameterized random graph model and assume that for any $\theta \in \Theta$, there exists $\mu_\theta$ the limiting ESD with respect to weak, almost sure convergence. We assume that $\mu_\theta\in \mathcal{M}_1(\Lambda)$, where $\Lambda$ is a bounded set. Consider $(\mathcal{G}_n)_{n\geq 1}$ a sequence of random graphs from distribution $P_{\theta^\star}$.  Let $\phi$ be a  kernel and $\sigma=\sigma_n$ a bandwidth that converges to 0. 
If the map $\theta \in \Theta \mapsto \mu_\theta(d\lambda) \in \mathcal{M}_1(\Lambda)$ is injective, continuous and $\Theta$ is compact, then the minimizer
\[\hat\theta_n = \argmin_{\theta \in \Theta} ||\mu^{\mathcal{G}_n}(\phi_{\cdot,\sigma}) - \mu_\theta(\phi_{\cdot,\sigma})||_1\]
converges in probability to $\theta^\star$ as $n \rightarrow \infty$. 
\end{theorem}

\begin{proof}
The proof is based on Theorem 3.2.8 from
\cite{DCD} about minimum contrast estimators. 
The first step is establishing the convergence of the contrast function 
$ D_n(\theta, \theta^\star) := \|\mu^{\mathcal{G}_n}(\phi_{\cdot,\sigma}) - \mu_{\theta}(\phi_{\cdot,\sigma})||_1$.  
Assuming that $\mu_\theta$ exists for any $\theta \in \Theta$, we get that for any $x \in \mathbb{R}$,  
\[ \mu^{\mathcal{G}_n}(\phi_{x,\sigma}) = \frac 1 n \sum_{i=1}^n \phi_{x,\sigma}(\lambda^{\mathcal{G}_n}_i/ \sqrt{n}) \mathop{\longrightarrow}_{n\to \infty} \mu_{\theta^\star}(\phi_{x,\sigma}), 
\]
$P_{\theta^\star}$-almost surely. Then, from dominated convergence theorem, we have 
\[
D_n(\theta, \theta^\star) 
\mathop{\longrightarrow}_{n\to \infty}  \| \mu_{\theta^\star}(\phi_{\cdot,\sigma}) - \mu_{\theta}(\phi_{\cdot,\sigma})||_1 := D(\theta, \theta^\star), 
\]
$P_{\theta^\star}$-almost surely.
Now, the second step is to establish that the limiting function $\theta \mapsto D(\theta, \theta^\star)$ has a strict minimum at $\theta^\star$.
Notice that the limit  $D(\theta,\theta^\star)$ is minimum if and only if $\mu_{\theta^\star}(\phi_{\cdot,\sigma})  =\mu_{\theta}(\phi_{\cdot,\sigma})$ almost surely (wrt the Lebesgue measure), namely 
\[ 
\forall x\in \mathbb{R}, \quad \int \phi[(x-\lambda)/\sigma] \mu_{\theta^\star}(d\lambda) = \int \phi[(x-\lambda)/\sigma] \mu_{\theta}(d\lambda)  .
\]
Classical properties of kernels  \citep[see for e.g. Proposition 1.2 in][]{tsybakov_book} imply that letting $\sigma \to 0$, we get 
\[ \mu_{\theta} (d\lambda)= \mu_{\theta^\star} (d\lambda).
\]
Since $\theta \mapsto \mu_\theta$ is injective, this implies $\theta=\theta^\star$. Thus, the limit $D(\theta,\theta^\star)$ is minimum if and only if $\theta=\theta^\star$. 

Finally, 
the map $\theta \mapsto  \mu_{\theta}(d\lambda)$ is continuous on the compact set $\Theta$ and thus uniformly continuous. This implies that \[
\forall \epsilon >0, \exists \eta >0  \text{ such that if } \|\theta-\theta'\| \leq \eta \text{ then }  |\mu_{\theta}(d\lambda) -  \mu_{\theta'}(d\lambda)| \leq \epsilon. 
\]
Thus we get that 
 \begin{align*}
| D_n(\theta,\theta^\star) -  D_n(\theta',\theta^\star)| &\leq \|\mu_\theta(\phi_{\cdot,\sigma}) - \mu_{\theta'}(\phi_{\cdot,\sigma})\|_1 \\
& \leq \sigma^{-1}\int \int_\Lambda \phi[(x-\lambda)/\sigma] \times  |\mu_\theta(d\lambda) -  \mu_{\theta'}(d\lambda)| dx\\
& \leq \int \int \phi(y) \times  |\mu_\theta(x-\sigma y) -  \mu_{\theta'}(x-\sigma y)| dy dx\\
& \leq \epsilon |\Lambda| ,
\end{align*}
where $|\Lambda| $ denotes the length of the interval $\Lambda$. 
For clarity, we assumed in the last line that $\mu_\theta$ is a measure with a density, but this has no consequence
as soon as $\|\theta-\theta'\| \leq \eta$. This establishes a uniform convergence (wrt $\theta$) of the contrast function $D_n$ as $n$ tends to infinity. Relying on Theorem 3.2.8 from \cite{DCD}, we obtain the desired result. 
\end{proof}

We can establish a similar result whether we consider cumulative distribution functions (CDF) of eigenvalues instead of spectral densities.

\begin{theorem}
\label{theo:cdf}
Let $\{ P_\theta ; \theta\in \Theta\}$ denote a parameterized random graph model and assume that for any $\theta \in \Theta$, there exists $\mu_\theta$ the limiting ESD with respect to weak almost sure convergence. We assume that $\mu_\theta\in \mathcal{M}_1(\Lambda)$, where $\Lambda$ is a bounded set. Consider $(\mathcal{G}_n)_{n\geq 1}$ a sequence of random graphs from distribution $P_{\theta^\star}$. 
If the map $\theta \in \Theta \mapsto \mu_\theta(d\lambda) \in \mathcal{M}_1(\Lambda)$ is injective, continuous and $\Theta$ is compact, then the minimizer
\[\hat\theta_n = \argmin_\theta ||F^{\mathcal{G}_n} - F_\theta||_1\]
converges in probability to $\theta^\star$ as $n \rightarrow \infty$.
\end{theorem}

\begin{proof}
The proof follows the same lines as the previous one, by first  noting that 
for any $x$, the function $\lambda \mapsto \1\{\lambda \le x\} $ is continuous almost everywhere. Thus  
\[ F^{\mathcal{G}_n}(x) \mathop{\longrightarrow}_{n\to \infty} F_\theta(x), 
\]
$P_{\theta^\star}$-almost surely. 
For the second step, notice that  $F_{\theta}$ entirely characterizes the distribution $\mu_\theta$. 
Finally, we similarly prove that 
 \begin{align*}
| D_n(\theta,\theta^\star) -  D_n(\theta',\theta^\star)| &\leq \|F_\theta - F_{\theta'}\|_1 \\
& \leq \int_{\Lambda} \int_{\Lambda} \1\{\lambda \le x\} |\mu_\theta(d\lambda) -  \mu_{\theta'}(d\lambda)| dx\\
& \leq \epsilon |\Lambda|^2  
\end{align*}
as soon as $\|\theta-\theta'\| \leq \eta$. 
\end{proof}

\section{Simulations}
\label{sec:sim}

In what follows, we describe experiments to study the block model spectral density, illustrate Theorems~\ref{theo:esd} and \ref{theo:cdf}, and evaluate the model selection procedures' performance. We used the Gaussian kernel to estimate the ESD and the Silverman's criterion~\cite{silverman_density_1986} to choose the bandwidth $\sigma$ in all experiments. To generate random graphs, we used package \emph{igraph} for  \emph{R}~\cite{Csrdi2006TheIS}.

First, in Section~\ref{sec:sim_block_model}, we analyze the ESD of BM random graphs. Indeed, the results in \cite{avrachenkov_spectral_2015} are limited to a specific BM and we  empirically explore the validity of the convergence of the ESD in more general BM. 
Then, in Section~\ref{sec:sim_theo}, we analyze the parameter estimator's empirical behavior to confirm results given by Theorems~\ref{theo:esd} and \ref{theo:cdf} in scenarios in which a limiting spectral density exists and has a known analytical form.
In Section~\ref{sec:sim_empirical}, we further explore scenarios in which the limiting spectral density may not exist. 
The distribution of the random graph model's eigenvalues is estimated empirically, as described in Algorithm~\ref{alg:model_selection}. Finally, in Section~\ref{sec:sim_model_selection}, we show experiments for the model selection procedure.

\subsection{Block Model spectral density}
\label{sec:sim_block_model}

We generated BM graphs randomly and compared the ESD of the adjacency matrix $\mathbf{A}/\gamma(n)$ with the centered matrix $\mathbf{\Tilde{A}}$ given in~\eqref{eq:Asbm} and the theoretical density derived by~\cite{avrachenkov_spectral_2015} that are recalled in Proposition~\ref{prop:avra}. We considered five different scenarios, described as follows.

\begin{enumerate}
\item \textbf{Scenario 1.} Graph generated by the same model described in  \cite{avrachenkov_spectral_2015} with a small number of blocks ($M=3$). We set each block size to $K = 300$. The probability of connecting vertices from different groups is set to $p_0 = 0.2$ and the probabilities within the groups $(p_1, p_2, p_3)$  to $(0.8, 0.5, 0.6)$, respectively.
\item \textbf{Scenario 2.} Graph generated by the same model described in  \cite{avrachenkov_spectral_2015} with a larger number of groups ($M=10$). We set each block size to $K = 300$. We set the probability to connect vertices from different groups to $p_0 = 0.2$ and the probabilities within the groups $(p_1, p_2, p_3, p_4, p_5, p_6, p_7, p_8, p_9, p_{10})$ to $(0.8, 0.5, 0.6, 0.7, 0.4, 0.9, 0.55, 0.42, 0.38, 0.86)$, respectively.
\item \textbf{Scenario 3.} Generalization of the model described by \cite{avrachenkov_spectral_2015} for different block sizes. We considered $M=3$ blocks of sizes 100, 80, and 300, respectively. We set the probability of connecting vertices from different groups to $p_0 = 0.2$, and the probabilities within the groups $(p_1, p_2, p_3)$ to $(0.8, 0.5, 0.6)$, respectively.
\item \textbf{Scenario 4.} Generalization of the model described by \cite{avrachenkov_spectral_2015} for $p_0$ larger than the probabilities within the groups. We considered $M=3$ blocks of size $K = 300$. We set the probability of connecting vertices from different groups to $p_0 = 0.9$ and the probabilities within the groups $(p_1, p_2, p_3)$ to $(0.8, 0.5, 0.6)$, respectively.
\item \textbf{Scenario 5.} Generalization of the model described by \cite{avrachenkov_spectral_2015} with different probabilities of connecting vertices from different groups. We considered $M=3$ blocks of size $K = 300$. We set the probabilities of connecting vertices from different groups (previously $p_0$) to $p_{12} = 0.1$, $p_{13} = 0.2$, and $p_{23} = 0.05$. The probabilities within the groups $(p_1, p_2, p_3)$ are set to  $(0.8, 0.5, 0.6)$, respectively.
\end{enumerate}

\begin{figure}[H]
\centering
   \includegraphics[width=\columnwidth]{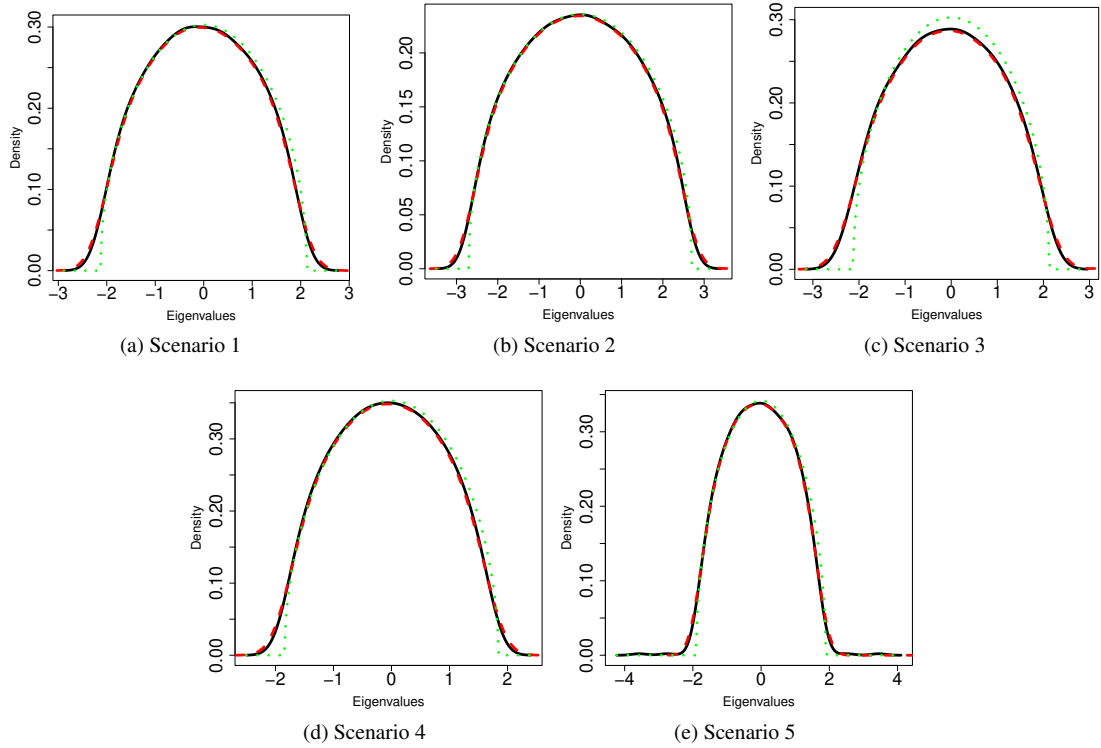}
     \caption{Comparison between the empirical spectral distribution (ESD) of the block model for matrices $\mathbf{\tilde{A}}$ (solid black line), $\mathbf{A}/\gamma(n)$ (red dashed line), and the theoretical distribution derived by \cite{avrachenkov_spectral_2015} (green dotted line). We considered different scenarios: (a) a small number of blocks, (b) a larger number of blocks, (c) blocks of different sizes, (d) probability between groups larger than the probabilities within the groups, and (e) different probabilities of connecting vertices from different groups. }
    \label{fig:esd_sbm}
\end{figure}

We compared the ESD of both the centered matrix $\mathbf{\tilde{A}}$ and the non-centered matrix $\mathbf{A}/\gamma(n)$ to \cite{avrachenkov_spectral_2015}'s theoretical density in Figure~\ref{fig:esd_sbm}. Note that the normalizing term $\gamma(n)$ depends on the unknown parameter value $p^\star$.
The density is computed numerically by inverting the Stieltjes transform.

For all the scenarios considered here, the ESDs of $\mathbf{A} /\gamma(n)$ and $\mathbf{\tilde{A}}$ are very close, suggesting that their ESD converges to the same distribution. 
For scenarios 1, 2, 4, and 5, the ESDs of $\mathbf{\tilde{A}}$ and $\mathbf{A}/\gamma(n)$ are close to the theoretical distribution derived by \cite{avrachenkov_spectral_2015}, suggesting that Avrachenkov \emph{et al.}'s results  could be further investigated for models with different probabilities to connect vertices from different groups, and large $p_0$. For scenario 3, the empirical densities are farther from the theoretical density. This might be due to numerical errors during ESD estimation or might indicate that Avrachenkov \emph{et al.}'s results do not apply to graphs with blocks of different sizes.

\subsection{Illustration of our convergence theorems for ER and DR models}
\label{sec:sim_theo}

To illustrate Theorems~\ref{theo:esd} and \ref{theo:cdf}, we considered scenarios in which we know the ESD analytic form, such as the ER and $d$-regular random graph models. 
Note however that in the latter case, the convergence is known for the adjacency matrix with no scaling, or scaled  with $\sqrt{d}$ while our theorems rely on a $\sqrt{n}$ normalization.
All experiments described in this section rely on the \emph{optimize} function from the \emph{R} package \emph{stats} to solve~\eqref{eq:parameter_estimator}, with precision $\epsilon=10^{-8}$. It uses a combination of golden section search and successive parabolic interpolation. 

\begin{figure}[H]
\centering
 \includegraphics[width=\textwidth]{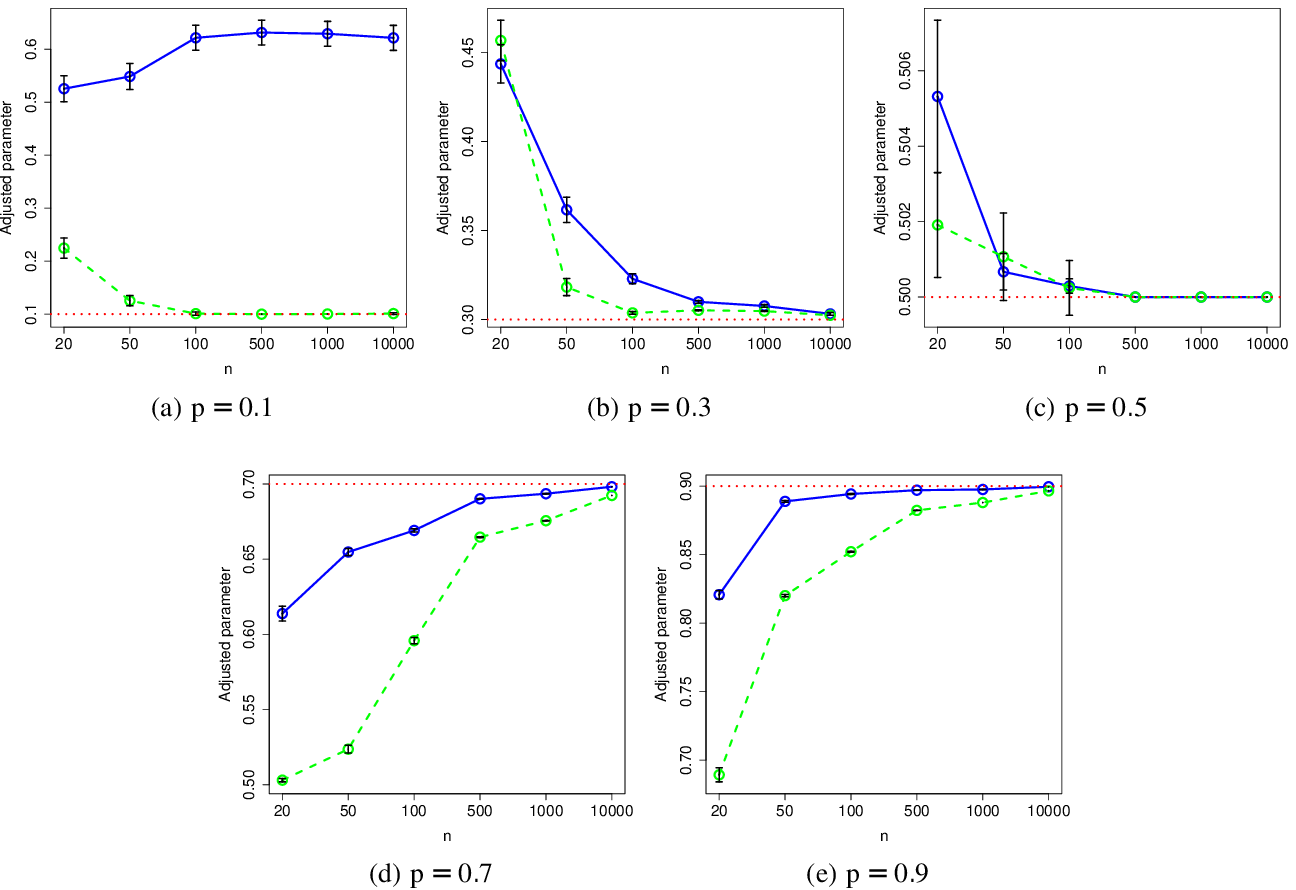}
 \caption{Average estimated parameter for the ER model. We generated $1\,000$ Erd\H{o}s-R\'{e}nyi random graphs of sizes $n = 20, 50, 100, 500, 1\,000, 10\,000$, and probability of connecting two vertices $p = 0.1$ (a), $p = 0.3$ (b), $p = 0.5$ (c), $p = 0.7$ (d), and $p = 0.9$ (e). In each plot, the $x$-axis and the $y$-axis correspond, respectively, to the graph size ($n$) and estimated parameter. The points and error bars indicate the average estimated parameters and the corresponding 95\% confidence intervals. The blue solid lines correspond to the results based on the $\ell_1$ distance between the ESD and the limiting theoretical spectral density of the ER model (semicircle law). The green dashed lines correspond to the results based on the eigenvalues' cumulative distribution. The red dotted lines indicate the true value.}
 \label{fig:sim_theo_er}
\end{figure}

For the ER model, we expect the ESD of $\mathbf{A}/\sqrt{n}$  to converge to the semicircle law when $p = \omega(n^{-1})$. In that case, the assumptions of Theorems~\ref{theo:esd} and~\ref{theo:cdf} hold when $\Theta = [0, 0.5]$ or $\Theta = [0.5, 1]$, and then the $\ell_1$ distance minimizer for problem~\eqref{eq:parameter_estimator} converges in probability to the true parameter. We generated $1\,000$ graphs of sizes $n = 20, 50, 100, 500, 1\,000, 10\,000$ with parameters $p = 0.1, 0.3, 0.5, 0.7, 0.9$. 
For $p < 0.5$, $p = 0.5$, and $p > 0.5$, we considered the intervals $\Theta = [0, 0.5]$, $\Theta = [0, 1]$, and $\Theta = [0.5, 1]$, respectively.  
In Figure~\ref{fig:sim_theo_er}, we show the average estimated parameter with a 95\% confidence interval for the ER model by using the spectral density (blue solid line) and the cumulative distribution (green dashed line). The red dotted line indicates the true parameter.

\begin{figure}[H]
\centering
 \includegraphics[width=\textwidth]{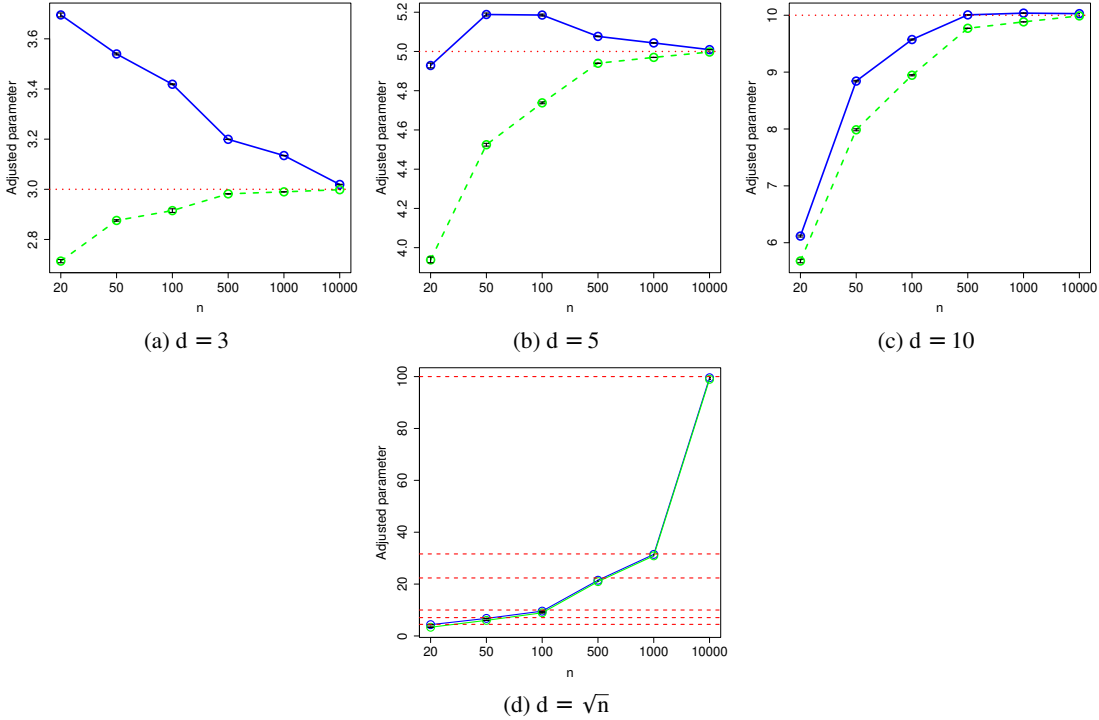}
 \caption{Average estimated parameter for the DR model. We generated $1\,000$  $d$-regular random graphs of sizes $n = 20, 50, 100, 500, 1\,000, 10\,000$, and degree $d = 3$ (a), $d = 5$ (b), $d = 10$ (c), and $d = \sqrt{n}$ (d). In each plot, the $x$-axis and the $y$-axis correspond, respectively, to the graph size ($n$) and estimated parameter. The points and error bars indicate the average estimated parameters and the corresponding 95\% confidence intervals. The blue solid lines correspond to the results based on the $\ell_1$ distance between the ESD and the limiting theoretical spectral density of the DR model (Kesten-McKay law). The green dashed lines correspond to the results based on the eigenvalues' cumulative distribution. The red dotted lines indicate the true value. }
 \label{fig:sim_theo_dr}
\end{figure}

Similarly, we performed simulation experiments with $d$-regular random graphs in scenarios in which the ESD of $\mathbf{A}$ 
converges weakly almost surely to the Kesten-McKay law. In what follows, we slightly modified the procedure from Algorithm~\ref{alg:model_selection} to rely on the spectrum of unnormalized adjacency matrices $\mathbf{A}$ (rather than $\mathbf{A}/\sqrt{n}$).
We generated $1\,000$ graphs of sizes $n = 20, 50, 100, 500, 1\,000, 10\,000$ with parameters $d = 3, 5, 10, \sqrt{n}$, namely fixed or slowly growing values $d$. Finally, we set $\Theta = [0, n]$.  Figure~\ref{fig:sim_theo_dr} shows the average estimated parameter with a 95\% confidence interval for the $d$-regular model by using the spectral density (blue solid line) and the cumulative distribution (green dashed line).

As expected by Theorems~\ref{theo:esd} and~\ref{theo:cdf}, we observe in Figures~\ref{fig:sim_theo_er} and~\ref{fig:sim_theo_dr} 
that the estimated parameter approximates the true parameter as the number of vertices increases.

\subsection{Convergence of the parameter estimates in other scenarios}
\label{sec:sim_empirical}
In this section, we explore the performance of the parameter estimation procedures based on $\ell_1$ distance (either through kernel or CDF). Our experiments include scenarios in which the assumptions of Theorems~\ref{theo:esd} and \ref{theo:cdf} may not hold, and a limiting spectral density does not exist, or it is unknown, such as the GRG, WS, and BA random graphs. In that case, to evaluate the model fitting procedure's performance based on the $\ell_1$ distance, we estimate the model's spectral density empirically, as described in Algorithm~\ref{alg:model_selection}. Our experiments also include BM graphs that \cite{takahashi_discriminating_2012} did not consider. For that model,  
we could compute the limiting spectral density numerically by inverting the Stieltjes transform (as is done in the previous section).
However, as it is computationally slow, we empirically estimate the BM spectral density for this experiment.

We generated $1\,000$ graphs of sizes $n = 50, 100, 500, 1\,000$ with parameters $r = 0.1, 0.3, 0.5, 0.7, 0.9$ for the GRG model, $p_r = 0.1, 0.3, 0.5, 0.7, 0.9$ for WS, and $p_s = 1.1, 1.3, 1.5, 1.7, 1.9$ for BA. For the two first models, in the search grid, we considered values between 0 and 1, varying with a step of size 0.001. For the BA model, we considered the interval $[1,4]$, and inspected values using a step size of 0.01. We considered two blocks of the same sizes and equal probabilities of connection within the groups for the BM graphs. We generated $1\,000$ graphs of sizes $n = 100, 600$.  We set the probability of connecting vertices from different groups to $p_0 = 0.2$ and the probabilities within the groups to $0.7$. Notice that for this model, the search space of the grid is bi-dimensional. We considered the space $[0.1,0.3]\times[0.6,0.8]$, with a step of $0.001$. For computing the eigenvalues, we used the $\sqrt{n}$ normalization of the adjacency matrix.

\begin{figure}[H]
\centering
\includegraphics[width=\textwidth]{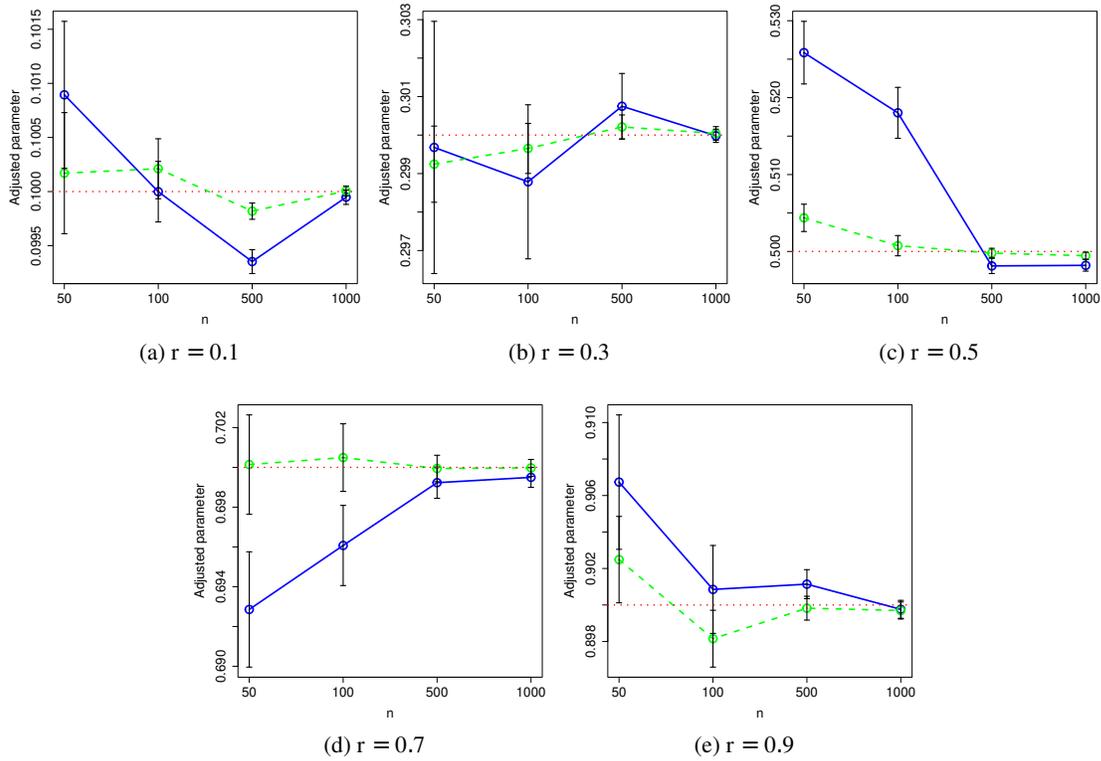}
 \caption{Average estimated parameter for the GRG model. We generated $1\,000$ GRG random graphs of sizes $n = 50, 100, 500, 1\,000$, and radius $r = 0.1$ (a), $r = 0.3$ (b), $r = 0.5$ (c), $r = 0.7$ (d), and $r = 0.9$ (e). 
In each plot, the $x$-axis and the $y$-axis correspond, respectively, to the graph size ($n$) and estimated parameter. The points and error bars indicate the average estimated parameters and the corresponding 95\% confidence intervals. The blue solid lines correspond to the results based on the $\ell_1$ distance between the ESD of the observed graph and the average ESD of the GRG model. The green dashed lines correspond to the results based on the eigenvalues' cumulative distribution. The red dotted lines indicate the true value.}
 \label{fig:sim_other_grg}
\end{figure}

\begin{figure}[H]
\centering
 \includegraphics[width=\textwidth]{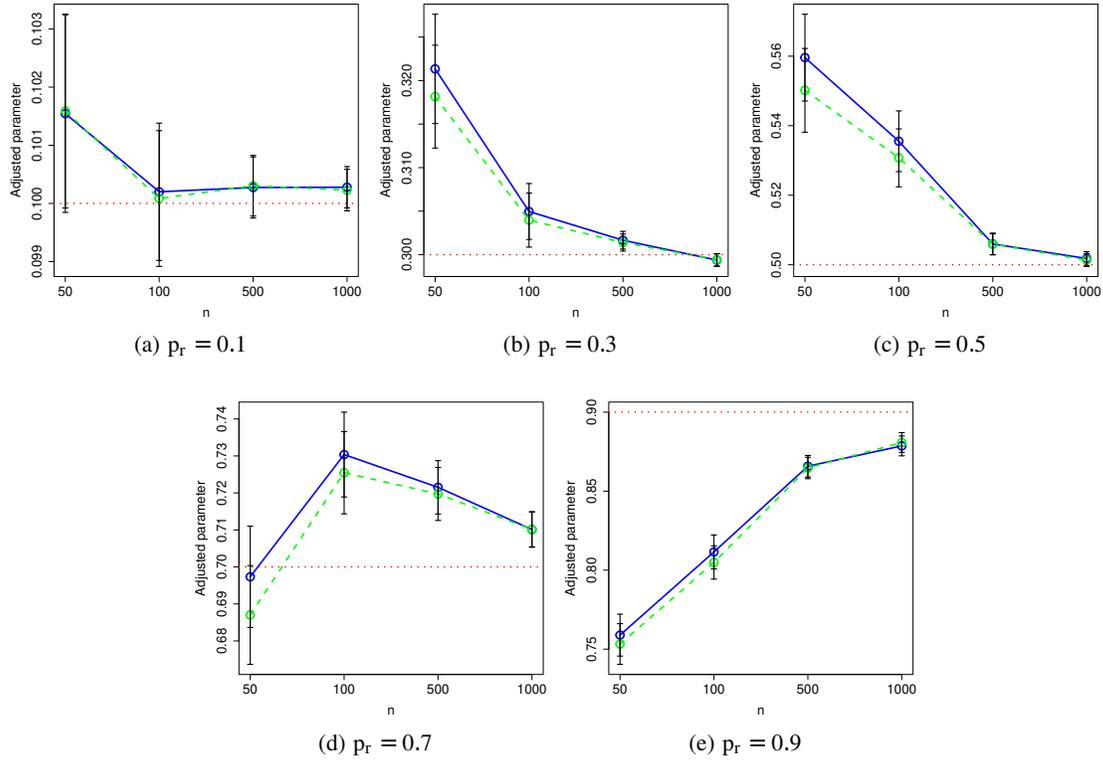}
 \caption{Average estimated parameter for the WS model. We generated $1\,000$ WS random graphs of sizes $n = 50, 100, 500, 1\,000$, and probability of reconnecting edges $p_r = 0.1$ (a), $p_r = 0.3$ (b), $p_r = 0.5$ (c), $p_r = 0.7$ (d), and $p_r = 0.9$ (e). In each plot, the $x$-axis and the $y$-axis correspond, respectively, to the graph size ($n$) and estimated parameter. The points and error bars indicate the average estimated parameters and the corresponding 95\% confidence intervals. The blue solid correspond to the results based on the $\ell_1$ distance  between the ESD of the observed graph and the average ESD of the WS model. The green dashed lines correspond to the results based on the eigenvalues' cumulative distribution. The red dotted lines indicate the true value.}
 \label{fig:sim_other_ws}
\end{figure}

\begin{figure}[H]
\centering
 \includegraphics[width=\textwidth]{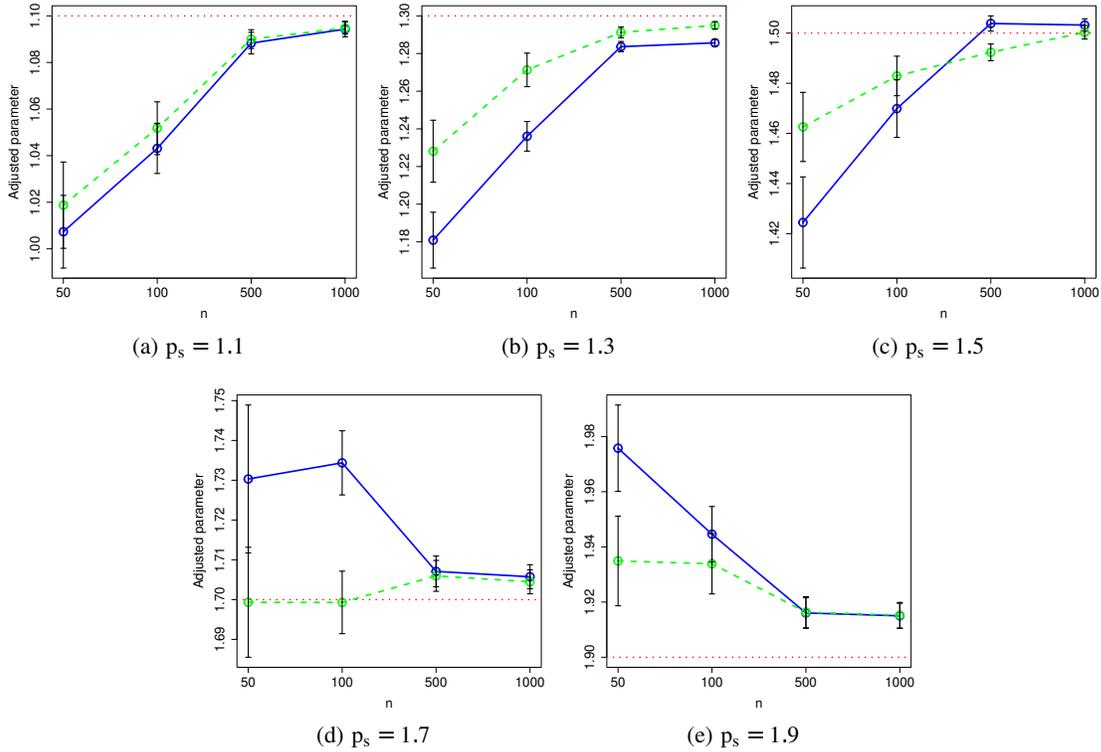}
 \caption{Average estimated parameter for the BA  model. We generated $1\,000$ BA random graphs of sizes $n = 50, 100, 500, 1\,000$, and  scale exponent  $p_s = 1.1$ (a), $p_s = 1.3$ (b), $p_s = 1.5$ (c), $p_s = 1.7$ (d), and $p_s = 1.9$ (e). In each plot, the $x$-axis and the $y$-axis correspond, respectively, to the graph size ($n$) and estimated parameter. The points and error bars indicate the average estimated parameters and the corresponding 95\% confidence intervals. The blue solid lines correspond to the results based on the $\ell_1$ distance between the ESD of the observed graph and the average ESD of the BA model. The green dashed lines correspond to the results based on the eigenvalues' cumulative distribution. The red dotted lines indicate the true value.}
 \label{fig:sim_other_ba}
\end{figure}

\begin{figure}[H]
\centering
 \includegraphics[width=0.7\textwidth]{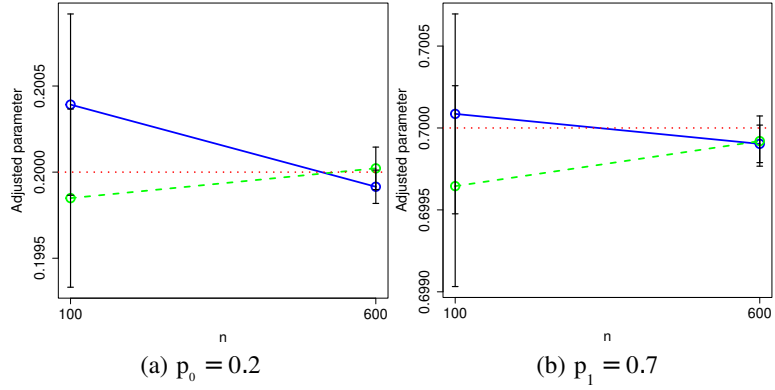}
 \caption{Average estimated parameter for the BM model with two blocks of the same sizes. We generated $1\,000$ BM random graphs of sizes $n = 100, 600$, with probability to connect vertices from different groups  $p_0 = 0.2$ (a), and probability within the clusters $p_1 = 0.7$ (b). In each plot, the $x$-axis and the $y$-axis correspond, respectively, to the graph size ($n$) and estimated parameter. The points and error bars indicate the average estimated parameters and the corresponding 95\% confidence intervals. The solid blue lines correspond to the results based on the $\ell_1$ distance between the ESD of the observed graph and the average ESD of the BM model. The green dashed lines correspond to the results based on the eigenvalues' cumulative distribution. The red dotted lines indicate the true value.}
 \label{fig:sim_other_bm}
\end{figure}

Figure~\ref{fig:sim_other_grg} shows the average estimated parameter and a 95\% confidence interval for the GRG model. Results obtained for the WS, BA, and BM models are shown in Figures~\ref{fig:sim_other_ws}, \ref{fig:sim_other_ba}, and \ref{fig:sim_other_bm} respectively. The blue solid lines correspond to the results based on the $\ell_1$ distance between the ESD of the observed graph and the average ESD of the model, estimated as described by Algorithm~\ref{alg:model_selection}. Similarly, the green dashed lines  correspond to the results based on the eigenvalues' empirical cumulative distribution. The red dotted lines indicate the true parameter.

In most scenarios, we observe that the estimated parameter approximates the actual value as the graph size increases. However, in some cases, results for $n=500$  and $n=1\,000$ are similar (see, e.g., BA model, when $p_s = 1.7, 1.9$). This limited convergence might be related to the fact that we only consider parameter values on the search grid (as described in Algorithm~\ref{alg:model_selection}).

\subsection{Model selection}
\label{sec:sim_model_selection}

To evaluate the performance of the model selection approach based on the $\ell_1$ distance between the ESD and its theoretical limit, we performed simulation experiments with $1\,000$ graphs of sizes $n = 50, 100, 500$ generated by the following models: Erd\H{o}s-R\'{e}nyi (ER) ($p = \beta$), Geometric (GRG) ($r = \beta$), $d$-regular (DR) ($d= 10\beta$), Watts-Strogatz (WS) ($p_r=\beta$), and Barab\'{a}si-Albert (BA) ($p_s = 1 + \beta$), where $\beta = 0.1, 0.5, 0.9$.  We also contrasted the performance against the original method by Takahashi et al.~\cite{takahashi_discriminating_2012}, based on the Kullback-Leibler divergence between spectral densities. Tables~\ref{tab:model_selection_0_1}, \ref{tab:model_selection_0_5}, \ref{tab:model_selection_0_9} show the confusion matrices obtained by the methods based on the spectral density and cumulative distribution of the eigenvalues for $\beta = 0.1, 0.5$, and $0.9$, respectively.

\begin{table}[H] 
\centering
\small
\caption{\textbf{Model selection for graphs generated with $\beta=0.1$}. The number of graphs correctly classified by the model selection approach based on the Kullback-Leibler divergence between empirical spectral densities (left)}, the $\ell_1$ distance between empirical spectral densities (middle), and the $\ell_1$ distance between empirical cumulative distributions (right). For each model, we generated $1\,000$ graphs of sizes $n = 50, 100, 500$. The set-up of the models was: the probability of connecting two vertices $p =\beta$ for the Erd\H{o}s-R\'{e}nyi (ER), the radius $r = \beta$ for the geometric (GRG),  the degree $d = 10\beta$  for the $d$-regular (DR), the probability of reconnecting edges $p_r =\beta$  for the Watts-Strogatz (WS), and the scale exponent $p_s = 1 + \beta$ for the Barab\'{a}si-Albert (BA).
\begin{tabular}
{>{\centering}m{0.01\textwidth}|>{\centering}m{0.05\textwidth}|>{\centering}m{0.12\textwidth}>{\centering}m{.12\textwidth}>{\centering}m{.12\textwidth}>{\centering}m{.12\textwidth}c}
\multirow{2}{*}{$n$} & \multirow{2}{*}{True} & \multicolumn{5}{c}{Selected model}\\
 & model &ER & DR & GRG & WS & BA\\
 \hline
& ER & \cellcolor{gray!25}$905|881|929$ & $6|7|5$ & $0|0|0$ & $89|112|66$ & $0|0|0$\\
& DR & $0|0|0$ & \cellcolor{gray!25}$1000|1000|1000$ & $0|0|0$ & $0|0|0$ & $0|0|0$\\
50 & GRG & $1|4|20$ & $0|0|0$ & \cellcolor{gray!25}$999|996|980$ & $0|0|0$ & $0|0|0$\\
& WS & $0|0|0$ & $0|0|0$ & $0|0|0$ & \cellcolor{gray!25}$1000|1000|1000$ & $0|0|0$\\
& BA & $39|108|16$ & $0|0|0$ & $0|3|0$ & $0|0|0$ & \cellcolor{gray!25}$961|889|984$\\
\hline
& ER & \cellcolor{gray!25}$1000|993|998$ & $0|7|2$ & $0|0|0$ & $0|0|0$ & $0|0|0$\\
& DR & $0|0|0$ & \cellcolor{gray!25}$1000|1000|1000$ & $0|0|0$ & $0|0|0$ & $0|0|0$\\
100 & GRG & $0|0|0$ & $0|0|0$ & \cellcolor{gray!25}$1000|1000|1000$ & $0|0|0$ & $0|0|0$\\
& WS & $0|0|0$ & $0|0|0$ & $0|0|0$ & \cellcolor{gray!25}$1000|1000|1000$ & $0|0|0$\\
& BA & $9|94|3$ & $0|0|0$ & $2|2|0$ & $0|0|0$ & \cellcolor{gray!25}$989|904|997$\\
\hline
& ER & \cellcolor{gray!25}$1000|994|998$ & $0|6|2$ & $0|0|0$ & $0|0|0$ & $0|0|0$\\
& DR & $0|0|0$ & \cellcolor{gray!25}$1000|1000|1000$ & $0|0|0$ & $0|0|0$ & $0|0|0$\\
500 & GRG & $0|0|0$ & $0|0|0$ & \cellcolor{gray!25}$1000|1000|1000$ & $0|0|0$ & $0|0|0$\\
& WS & $0|0|0$ & $0|0|0$ & $0|0|0$ & \cellcolor{gray!25}$1000|1000|1000$ & $0|0|0$\\
& BA & $0|2|0$ & $0|0|0$ & $0|0|0$ & $0|0|0$ & \cellcolor{gray!25}$998|1000$\\
\hline
\end{tabular}\hspace*{-1.5cm}
\label{tab:model_selection_0_1}
\end{table}

\begin{table}[H] 
\centering
\small
\caption{\textbf{Model selection for graphs generated with $\beta=0.5$.} The number of graphs correctly classified by the model selection approach based on the Kullback-Leibler divergence between empirical spectral densities (left)}, the $\ell_1$ distance between empirical spectral densities (middle), and the $\ell_1$ distance between empirical cumulative distributions (right). For each model, we generated $1\,000$ graphs of sizes $n = 50, 100, 500$. The set-up of the models was: the probability of connecting two vertices $p =\beta$ for the Erd\H{o}s-R\'{e}nyi (ER), the radius $r = \beta$ for the geometric (GRG),  the degree $d = 10\beta$  for the $d$-regular (DR), the probability of reconnecting edges $p_r =\beta$  for the Watts-Strogatz (WS), and the scale exponent $p_s = 1 + \beta$ for the Barab\'{a}si-Albert (BA).
\begin{tabular}
{>{\centering}m{0.01\textwidth}|>{\centering}m{0.05\textwidth}|>{\centering}m{0.12\textwidth}>{\centering}m{.12\textwidth}>{\centering}m{.12\textwidth}>{\centering}m{.12\textwidth}c}
\multirow{2}{*}{$n$} & \multirow{2}{*}{True} & \multicolumn{5}{c}{Selected model}\\
 & model &ER & DR & GRG & WS & BA\\
 \hline
& ER & \cellcolor{gray!25}$947|948|911$ & $53|52|89$ & $0|0|0$ & $0|0|0$ & $0|0|0$\\
& DR & $0|1|0$ & \cellcolor{gray!25}$1000|999|1000$ & $0|0|0$ & $0|0|0$ & $0|0|0$\\
50 & GRG & $0|2|0$ & $0|0|0$ & \cellcolor{gray!25}$1000|998|1000$ & $0|0|0$ & $0|0|0$\\
& WS & $0|21|0$ & $0|0|0$ & $0|0|0$ & \cellcolor{gray!25}$1000|979|1000$ & $0|0|0$\\
& BA & $2|64|3$ & $0|0|0$ & $3|103|3$ & $0|0|0$ & \cellcolor{gray!25}$995|833|994$\\
\hline
& ER & \cellcolor{gray!25}$966|949|951$ & $53|51|49$ & $0|0|0$ & $0|0|0$ & $0|0|0$\\
& DR & $0|0|0$ & \cellcolor{gray!25}$1000|1000|1000$ & $0|0|0$ & $0|0|0$ & $0|0|0$\\
100 & GRG & $0|0|0$ & $0|0|0$ & \cellcolor{gray!25}$1000|1000|1000$ & $0|0|0$ & $0|0|0$\\
& WS & $3|6|0$ & $0|0|0$ & $0|0|0$ & \cellcolor{gray!25}$997|994|1000$ & $0|0|0$\\
& BA & $1|121|1$ & $0|0|0$ & $95|125|1$ & $0|0|0$ & \cellcolor{gray!25}$904|754|998$\\
\hline
& ER & \cellcolor{gray!25}$899|908|922$ & $101|92|78$ & $0|0|0$ & $0|0|0$ & $0|0|0$\\
& DR & $0|0|0$ & \cellcolor{gray!25}$1000|1000|1000$ & $0|0|0$ & $0|0|0$ & $0|0|0$\\
500 & GRG & $0|0|0$ & $0|0|0$ & \cellcolor{gray!25}$1000|1000|1000$ & $0|0|0$ & $0|0|0$\\
& WS & $0|0|0$ & $0|0|0$ & $0|0|0$ & \cellcolor{gray!25}$1000|1000|1000$ & $0|0|0$\\
& BA & $0|31|0$ & $0|0|0$ & $0|11|0$ & $0|0|0$ & \cellcolor{gray!25}$1000|958|1000$\\
\hline
\end{tabular}\hspace*{-1.5cm}
\label{tab:model_selection_0_5}
\end{table}

\begin{table}[H] 
\centering
\small
\caption{\textbf{Model selection for graphs generated with $\beta=0.9$.} The number of graphs correctly classified by the model selection approach based on \blue{the Kullback-Leibler divergence between empirical spectral densities (left)}, the $\ell_1$ distance between empirical spectral densities (middle), and the $\ell_1$ distance between empirical cumulative distributions (right).For each model, we generated $1\,000$ graphs of sizes $n = 50, 100, 500$. The set-up of the models was: the probability of connecting two vertices $p =\beta$ for the Erd\H{o}s-R\'{e}nyi (ER), the radius $r = \beta$ for the geometric (GRG),  the degree $d = 10\beta$  for the $d$-regular (DR), the probability of reconnecting edges $p_r =\beta$  for the Watts-Strogatz (WS), and the scale exponent $p_s = 1 + \beta$ for the Barab\'{a}si-Albert (BA).} 
\begin{tabular}
{>{\centering}m{0.01\textwidth}|>{\centering}m{0.05\textwidth}|>{\centering}m{0.12\textwidth}>{\centering}m{.12\textwidth}>{\centering}m{.12\textwidth}>{\centering}m{.12\textwidth}c}
\multirow{2}{*}{$n$} & \multirow{2}{*}{True} & \multicolumn{5}{c}{Selected model}\\
 & model &ER & DR & GRG & WS & BA\\
 \hline
& ER & \cellcolor{gray!25}$1000|1000|1000$ & $0|0|0$ & $0|0|0$ & $0|0|0$ & $0|0|0$\\
& DR & $0|56|4$ & \cellcolor{gray!25}$1000|944|996$ & $0|0|0$ & $0|0|0$ & $0|0|0$\\
50 & GRG & $0|71|0$ & $0|0|0$ & \cellcolor{gray!25}$1000|929|1000$ & $0|0|0$ & $0|0|0$\\
& WS & $19|133|22$ & $0|0|0$ & $0|0|0$ & \cellcolor{gray!25}$981|867|978$ & $0|0|0$\\
& BA & $0|9|0$ & $0|0|0$ & $3|22|0$ & $0|0|0$ & \cellcolor{gray!25}$997|969|1000$\\
\hline
& ER & \cellcolor{gray!25}$1000|1000|1000$ & $0|0|0$ & $0|0|0$ & $0|0|0$ & $0|0|0$\\
& DR & $0|1|0$ & \cellcolor{gray!25}$1000|999|1000$ & $0|0|0$ & $0|0|0$ & $0|0|0$\\
100 & GRG & $0|0|0$ & $0|0|0$ & \cellcolor{gray!25}$1000|1000|1000$ & $0|0|0$ & $0|0|0$\\
& WS & $5|124|59$ & $0|0|0$ & $0|0|0$ & \cellcolor{gray!25}$995|876|941$ & $0|0|0$\\
& BA & $0|0|0$ & $0|0|0$ & $2|2|0$ & $0|0|0$ & \cellcolor{gray!25}$998|998|1000$\\
\hline
& ER & \cellcolor{gray!25}$1000|1000|1000$ & $0|0|0$ & $0|0|0$ & $0|0|0$ & $0|0|0$\\
& DR & $0|0|0$ & \cellcolor{gray!25}$1000|1000|1000$ & $0|0|0$ & $0|0|0$ & $0|0|0$\\
500 & GRG & $0|0|0$ & $0|0|0$ & \cellcolor{gray!25}$1000|1000|1000$ & $0|0|0$ & $0|0|0$\\
& WS & $0|0|0$ & $0|0|0$ & $0|0|0$ & \cellcolor{gray!25}$1000|1000|1000$ & $0|0|0$\\
& BA & $0|0|0$ & $0|0|0$ & $0|0|0$ & $0|0|0$ & \cellcolor{gray!25}$1000|1000|1000$\\
\hline
\end{tabular}\hspace*{-1.5cm}
\label{tab:model_selection_0_9}
\end{table}

Tables~\ref{tab:model_selection_0_1},~\ref{tab:model_selection_0_5}, and~\ref{tab:model_selection_0_9} show that the performance is similar across the methods. The number of hits increases as the graphs become larger. The only exception is when $\beta = 0.5$ for ER graphs. In that case, the method confounds some ER graphs with DR graphs.
Note that the limiting spectral densities have the form of semicircle laws (though not obtained through the same scaling of the adjacency matrix).

\section{Real world networks}
\label{sec:real_world_networks}

To illustrate applications of the model selection approach based on the $\ell_1$ distance between spectral densities, we applied the procedure to three real-world networks, described as follows.

\begin{enumerate}
    \item Yeast protein-protein interaction (PPI) network. It contains $1\,966$ vertices corresponding to proteins and $2\,705$ edges indicating whether two proteins physically interact within the cell. We downloaded data from \url{http://networksciencebook.com/} \cite{yu_high-quality_2008}.
    \item E. coli bacteria metabolic network. It contains $1\,039$ vertices corresponding to metabolites and $5\,802$ edges indicating reactions between pairs of metabolites. We downloaded data from \url{http://networksciencebook.com/} \cite{schellenberger_bigg_2010}.
    \item Facebook ego network. From 10 Facebook ego networks (each one containing friends of a starting user), we selected the one with more than 1\,000 vertices. The selected network has $1\,034$ vertices corresponding to users and $53\,498$ edges indicating friendships on Facebook. We downloaded the data from \url{https://snap.stanford.edu/data/ego-Facebook.html} \cite{mcauley_learning_2012}.
\end{enumerate}

We considered a list of four candidate models: ER (an utterly random procedure),  GRG, WS, and BA (models with few parameters that share characteristics of complex networks). Table~\ref{tab:real_world_networks} shows the  $\ell_1$ distance between the ESD of the actual network and the average ESD of each model.

\begin{table}[H]
  \centering
  \caption{Model selection for real-world networks. We applied the model selection procedure for three different real networks (yeast protein-protein interaction network, E. coli metabolic network, and Facebook ego network). We considered four different candidate models: Erd\H{o}s-R\'enyi (ER), geometric random graph (GRG), Watts-Strogatz (WS), and Barab\'asi-Albert (BA). The values indicate the $\ell_1$ distance between the empirical spectral density of the network and the model. We highlighted in bold the smallest value for each network.}
  \begin{tabular}{c|cccc}
    Network & ER & GRG & WS & BA\\ 
    \hline
    Yeast PPI & 0.941 & 0.433 & 0.486 & \textbf{0.108}\\
   E. coli metabolic network & 0.636 & 0.367 & 0.455 & \textbf{0.17}\\
   Facebook network & 0.296 & \textbf{0.207} & 0.238 & 0.221\\
    \hline
  \end{tabular}
   \label{tab:real_world_networks}
\end{table}

Notice that BA is the only model from our list that builds scale-free graphs, a prevalent pattern in biological networks. For example, many protein-protein interaction networks and metabolic networks are known to share properties with scale-free graphs~\cite{jeong_lethality_2001, jeong_large-scale_2000}. Indeed, for both the Yeast PPI and the E. coli metabolic networks, BA is the model with the most similar spectral density. Similarly, the original procedure by \cite{takahashi_discriminating_2012}, based on the Kullback-Leibler divergence between spectral densities, also selects BA (among ER, WS, and BA) for PPI networks of different organisms.

Among the three networks, the Facebook ego network is the least similar to scale-free graphs. It contains social circles corresponding to densely connected groups of vertices~\cite{mcauley_learning_2012}. When clusters are well represented spatially, the graph may share characteristics of a geometric graph. In our analysis, the model selection approach selected GRG for this Facebook ego graph.

\section{Discussion}
\label{sec:discussions}

By reviewing results on the ESD of random graphs' convergence properties, we observed that the Wigner semicircle law could be applied to the ER and DR models when the adjacency matrices are adequately normalized, as shown in Table~\ref{tab:esd_convergence}. For the BM random graph, we can express the limiting ESD of the centered and normalized matrix through the Stieltjes transform under certain conditions (blocks with the same size, same probabilities of connecting vertices from different groups, and inter-groups probability smaller than the  intra-groups probabilities). Besides, our simulation experiments suggested that convergence of the ESD of non-centered adjacency matrices and centered BM adjacency matrices are similar. The ESD may also converge to the limiting ESD obtained through the Stieltjes transform when the probabilities to connect vertices from different groups are different, and the inter-groups probability  $p_0$ is larger than the intra-groups probabilities $p_m$.

\begin{table}[H]
  \centering
    \caption{Convergences of ESD in different models. }
  \begin{tabular}{ccc}
    Model & Matrix & Limit\\ 
    \hline
   ER & $\mathbf{A}/\sqrt{n}$ 
   & $\mu_{p, sc}(x)
=\frac 1 {2 \pi p(1-p)} \sqrt{[4p(1-p)-x^2]_{+}} $\\
& $\mathbf{A}/\sqrt{np(1-p)}$ 
& $\mu_{sc}(x) =
 \frac 1 {2  \pi}  \sqrt{[4-x^2]_{+}} $\\
    \hline
    DR & $\mathbf{A}$
    & $\mu_{d}(x) = \frac{d}{{2\pi (d^2 -x^2)}} \sqrt{[4(d-1)-x^2]_+}  $\\
& $\mathbf{A}/\sqrt{d-1}$ 
&   $\tilde \mu_{d}(x) =\left( 1 + \frac 1 {d-1} - \frac {x^2} {d} \right)^{-1} \frac{1}{2\pi}\sqrt{[4-x^2]_{+}} $ \\
    \hline
    BM &$(\mathbf{A}-\mathbb{E}(\mathbf{A}))/\sqrt{np^\star(1-p^\star)}$&  \text{expression through Stieltjes transform} \\
    \hline
  \end{tabular}
   \label{tab:esd_convergence}
\end{table}

The ESDs of the geometric random graph and the linear preferential attachment models converge weakly almost surely to a limiting density function as $n \to \infty$. However, an analytical expression for this limiting function is unknown. 

There is still no theoretical result on the ESD convergence for other models, such as the WS and BA. However, simulation experiments suggest that the ESD of those models approximate to a fixed distribution as the number of vertices increases \cite{farkas_spectra_2001}.

Therefore, there are scenarios in which we may assume that the ESD converges weakly to a limiting distribution, as required by Theorems~\ref{theo:esd} and~\ref{theo:cdf}. If all assumptions hold, then the minimizer of the $\ell_1$ distance between the observed graph's eigenvalue distribution and the limiting ESD of the model converges in probability to the true parameter as $n \to \infty$.

When assumptions of Theorems~\ref{theo:esd} and~\ref{theo:cdf} are satisfied, our simulations in Section~\ref{sec:sim_theo} show that the estimated parameters approximate the true values as the graph size increases. When we can not obtain a limiting spectral distribution, we estimate it empirically by randomly sampling graphs from the random graph model (as described in Algorithm~\ref{alg:model_selection}). In this case, the estimated parameter also approximates the true parameter as the number of vertices increases, as shown by simulation experiments in Section~\ref{sec:sim_empirical}. 

The model-fitting procedure's performance depends on the choice of the ESD or the CDF. The convergence of the estimated parameter based on the ESD, in general, is faster for the ER model. The model-fitting procedures' performance based on the ESD and CDF was similar for the WS model. For the GRG and BA random graph models, we observed a better performance by using the CDF. The reason is that sometimes ESD estimation/bandwidth choice is difficult. Thus, in that case, using the cumulative distribution of the eigenvalues instead of the spectral density is more convenient.

One of the main limitations of the procedure in Algorithm~\ref{alg:model_selection} is its high computational cost. 
In the following complexity analysis, we neglect the number of models $N$ that is supposed to be fixed (and small in general).

The spectral decomposition has a complexity $O(n^3)$. Thus, when the analytical form of the limiting spectral distribution is known, the complexity of Algorithm~\ref{alg:model_selection} is $O(n^3)$, while when we use a Monte Carlo approximation for the limiting ESD, it increases to  $O(n^3|\tilde\Theta| M)$. Future research for reducing computational cost of Algorithm~\ref{alg:model_selection} could explore either a faster approximation of the spectrum of each  graph and thus its ESD, or an approximation of this ESD. Indeed, a faster approximation of the ESD does not require the direct computation of all eigenvalues~\cite{newman_spectra_2019, Lin_etal_2016}, and has complexity of the order of $O(|E|C)$, where $C$ is associated to the number of random vectors and number of moments used by the approximation algorithm. For large graphs, we  have $C << |E|$. Then, if this approximation is used in our model fitting procedure, the final computational cost could be reduced to $O(|E|)$ (when analytic form of limiting ESD is known) or $O(|E||\tilde\Theta| M)$ (when MC approximation is needed).
Notice, however, that more experiments are necessary to validate these approximations.

Furthermore, the procedure in Algorithm~\ref{alg:model_selection} might not be the most recommended approach for classical models such as the ER and BM. The reason is that there are already accurate and efficient parameter estimators \cite{snijders_estimation_1997, ambroise_new_2012} for them. However,  Algorithm~\ref{alg:model_selection} gains in generality. Every time we fix a parametric model, we can build a procedure to estimate the parameter and a divergence between the model and the observed graph. This generality is particularly useful when considering a list of random graph models, and we want to choose the one that best approximates the observed graph.

Indeed, the experiments in Section~\ref{sec:sim_model_selection} suggest that choosing the model that minimizes the $\ell_1$ distance between spectral densities may present a hit rate of 100\% when the graph size is large enough. However, it is essential to notice that certain combinations of parameters may generate similar graphs from different models, such as the ER and DR graphs.

The illustration of the model selection approach on real-world networks (Section~\ref{sec:real_world_networks}) suggested that results are consistent with known properties of biological and social networks. However, it is essential to notice that this procedure does not correspond to a goodness-of-fit test. It might happen that the best model does not fit that well the observed network and is only the least bad. 
It would be interesting to test whether a given model is likely to generate the observed spectral density in future work. 

We focused on studying the ESD of the graph adjacency matrix. However, assumptions of Theorems~\ref{theo:esd} and \ref{theo:cdf} may also hold for other graph matrices, such as the Laplacian and normalized Laplacian, as discussed by \cite{ding_spectral_2010, gu_spectral_2016}. Therefore, results are general in terms of different random graph models and different graph matrices whose ESDs satisfy the assumptions of  Theorems~\ref{theo:esd} and \ref{theo:cdf}.

\section{Conclusions}
\label{sec:conclusions}

We reviewed the convergence properties of random graph's ESDs. Our review and simulation study suggest that different random graph models' ESD converge to a limiting density function. Based on the ESD's convergence properties, we prove that the minimizer of the $\ell_1$ distance of the ESD and the limiting density function is a consistent estimator of the random graph model parameter. This result also holds for the CDF of the eigenvalues. The main advantage of this model-fitting approach is its generality.

Furthermore, our simulation experiments  suggest that the procedures proposed by Takahashi et al. \cite{takahashi_discriminating_2012} may be modified to rely on $\ell_1$ distance between eigenvalue distribution functions to select a random graph model. 

In our computational experiments we used code from the \emph{statGraph} package available at CRAN (\url{https://cran.r-project.org/web/packages/statGraph/}). A zip file containing the R codes of all the experiments are available at \url{https://www.ime.usp.br/~suzana/code_spectral_density.zip}.

\section*{Funding}

This work was supported by S\~ao Paulo Research Foundation [grant numbers 2015/21162-4, 2017/12074-0, 2018/21934-5, 2019/22845-9, 2020/08343-8]; National Council for Scientific and Technological Development [grant number 303855/2019-3]; Coordination for the Improvement of Higher Education Personnel (Finance code 001); Alexander von Humboldt Foundation; and the Academy of Medical Sciences - Newton Fund.


\bibliographystyle{comnet}
\bibliography{references}
%


\end{document}